\documentclass[reqno, 11pt]{amsart}

\newtheorem{theorem}{Theorem}[section]
\newtheorem{lemma}[theorem]{Lemma}
\newtheorem{corollary}[theorem]{Corollary}

\theoremstyle{definition}
\newtheorem{assumption}[theorem]{Assumption}
\newtheorem{definition}[theorem]{Definition}
  \newtheorem{exercise*}[theorem]{Exercise*}

\theoremstyle{remark}

\newtheorem{remark}[theorem]{Remark}

\makeatletter
\def\dashintindex{\operatorname%
{-\kern-.7em\DOTSI\intop\ilimits@}}%
\def\dashint{\operatorname%
{\,\,\text{\bf--}\kern-.98em\DOTSI\intop\ilimits@\!\!}}
\makeatother

\newcommand\sfu{{\sf u}}

\newcommand\bB{\mathbb{B}}
\newcommand\bC{\mathbb{C}}

\newcommand\bM{\mathbb{M}}
\newcommand\bR{\mathbb{R}}
\newcommand\bQ{\mathbb{Q}}
\newcommand\bS{\mathbb{S}}
\newcommand\bZ{\mathbb{Z}}

 \newcommand{\loc}{{\rm loc}}

\newcommand\cF{\mathcal{F}}

\newcommand\cL{\mathcal{L}}
\newcommand\cM{\mathcal{M}}

\newcommand\cbrk{\text{$]$\kern-.15em$]$}}
\newcommand\opar{
\text{\,\raise.2ex\hbox{${\scriptstyle |}$}\kern-.34em$($}}

\newcommand{\mysection}[1]{\section{#1}
 \setcounter{equation}{0}}

\begin{document}

\title[Extrapolation theorem and applications]
{  Rubio de Francia extrapolation theorem
and related topics in the theory of
 elliptic and parabolic equations. A survey}
\author{N.V. Krylov}
\address{127 Vincent Hall, University of Minnesota, Minneapolis, MN, 55455}
\email{nkrylov@umn.edu}

\renewcommand{\subjclassname}{%
\textup{2000} Mathematics Subject Classification}

\date{}%
\subjclass{35-02, 35J99, 35K55}
\keywords{Rubio de Francia, extrapolation theorem,
mixed norm estimates}

\begin{abstract} 
We give a brief overview of the history 
of the Sobolev mixed norm theory of
linear   elliptic and 
parabolic equations and the recent development
in this theory
based on the Rubio de Francia extrapolation theorem.
A self contained proof of this theorem
along with other relevant tools
of Real Analysis are also presented
as well as an application to mixed norm
estimates for fully nonlinear equations.

\end{abstract}

\maketitle

\newcommand{\WO}{\overset{\scriptscriptstyle0}%
{W}\,\!}

\mysection{Introduction}
                                 \label{section 12.19.1}

The goal of this paper is to show how
the extrapolation theorem of Rubio de Francia
combined with more or less standard techniques
from the theory of partial differential equations
allows one to get estimates for solutions 
in Sobolev spaces with mixed norms ``for free''. 
We give a complete proof of this theorem in spite
of the fact that it can be found in a few articles
and books. The reason is that, if you try to learn
how this extrapolation theorem is proved by reading
books or articles on real or harmonic analysis, 
you will be buried under mountains of very beautiful
and fascinating results and it is not an easy task
to sort out which only very few of them are actually needed to
prove the extrapolation theorem.
We collected these few with complete proofs
including the extrapolation theorem on twelve pages.

The author's interest in equations in spaces with mixed
norms arose in connection with stochastic partial
differential equations where these norms and embedding
theorems allow one to obtain useful information
on solutions (see \cite{Kr_00}). 
See also \cite{GS_91} and \cite{MS_95}
for applications to the Navier-Stokes
equations   and \cite{SW_06} and \cite{We_02} for
applications to other problems.
The spaces used in \cite{Kr_00}  
were of type $L_{q}(\bR,H^{\gamma}_{p,\theta}(\bR^{d}_{+}))$,
$q\geq p$,
where $H^{\gamma}_{p,\theta}(\bR^{d}_{+})$ is an $H^{\gamma}_{p}$
type space with weights allowing certain blow up
near the boundary of the half space $\bR^{d}_{+}$.
In \cite{Kr_00} there is a restriction $q\geq p\geq2$
dictated by the presence of stochastic terms.
In \cite{Kr_01_1} stochastic terms from \cite{Kr_00}
are eliminated and $p,q$ are arbitrary in $(1,\infty)$.
The range of weights in \cite{Kr_01_1}, dictated
by the applications to stochastic partial
differential equations, was later
extended to the optimal one in \cite{KN_09}.

If there are no stochastic terms and
the coefficients are independent of time,
then we are dealing with the heat equation.
As far as we know,  the first explicitly stated
result on a priori estimates for the heat equation
 in spaces $L_{q}([0,T],W^{2}_{p}(\bR^{d}))$ with
$p,q\in(1,\infty)$
appeared in \cite{MS_95} with a somewhat sketchy proof
  based on the Calder\'on-Zygmund
theorem for operator-valued kernels which at that time
was only published for locally summable kernels
and formally was not applicable. Complete proofs
for more general case of the coefficients which depend
only on $t$ appeared in \cite{Kr_01_2} (one of the first
papers with mixed norms and the coefficients measurable in $t$).

Probably the earliest mentioning of
the results of that kind appeared in
\cite{So_64}. The proof in \cite{MS_95} uses
  the classical result for $q=p$ using which
  the above mentioned a priori estimate
from \cite{MS_95} can be easily obtained from a   general
result of \cite{CV_86} proved
(and stated earlier in \cite{So_64}) for abstract analytic
semigroups. This a priori estimate is also contained as a
particular case in Theorem 2.1 of \cite{GS_91} 
which is proved on the basis of a general theorem
on invertibility of the sum of two resolvent-commuting
operators
(see also \cite{GGS_91}, \cite{GGS_93}). 
 
Before that the case $p\in(1,\infty),q=1$ was
considered in \cite{Le_67}.
Several proofs and generalizations of the result from \cite{MS_95}
were given in \cite{Kr_01_1} where equations 
with measurable in time coefficients in $\bR^{d}$
and $\bR^{d}_{+}$ are treated and in the latter case
again weights are introduced.

In \cite{Kr_01_2}
a systematic approach to parabolic equations in $\bR^{d+1}$
in spaces $L_{q}(\bR,W^{2}_{p}(\bR^{d}))$ is given on the basis
of the Calder\'on-Zygmund theorem.
The coefficients are supposed to be independent of the space variable
and measurable with respect to the time variable.
  If $p=\infty$,
it is natural to replace the space
$L_{q}(\bR,W^{2}_{p}(\bR^{d}))$ with
$L_{q}(\bR,C^{2+\alpha}(\bR^{d}))$ and ask
if the results still hold.
In \cite{Kr_02} we prove that the answer
  to this question is indeed positive.
The approach in \cite{Kr_02} is quite different
from \cite{CV_86}  
and is based on simple  estimates of the heat potentials
and well-known properties
of the Hardy-Littlewood maximal function. 

  It is also worth mentioning the articles
\cite{CD_00} and \cite{HP_97}
containing a quite extensive references
where similar issues
are treated
in a very general setting of homogeneous spaces
and analytic semigroups.  
The semigroup approach or the approach based
on resolvent-commuting operators seem to produce
no results for parabolic equations
with coefficients measurable in time.
Good sources of references
and discussions of methods and obtained results
are also found in
  \cite{DHP_03}, \cite{HHH_06},    
  \cite{SW_06}, and \cite{We_02}.

Up until recently, excluding
\cite{Kr_02}, \cite{Kr_07_1},
  \cite{SW_06},  and \cite{We_02} and a few references therein,
in most other papers concerning $L_{q}(L_{p})$-spaces
the methods heavily depend on the properties
of the elliptic part of the operators,
which are supposed to be independent of
$t$ and have well behaving resolvent
or generate a ``good" semigroup.
However, in \cite{Am_04} (also see references therein)
 there is a general
theorem allowing one to treat the case when
the coefficients are continuous in $t$.
These restrictions exclude
parabolic equations with coefficients measurable or even
VMO in $t$
(even if they are independent of $x$,
the case considered in \cite{Kr_02}).
In particular, in \cite{HHH_06} the authors only consider
equations with VMO coefficients independent
of time, although combining
their results with \cite{Am_04} would include equations with
coefficients continuous in $t$. 
By  the way, in the particular case
that $q=p$ this also does not allow one to cover
the results of \cite{BC_93}, where
the coefficients  are in $VMO(\bR^{d+1})$.
Speaking about the case $q=p$, which is not
our main interest here, it is worth saying
that there is a quite extensive literature
about linear equations and systems with VMO coefficients.
The interested reader can consult  \cite{BW_10}, \cite{BPW_13},
 \cite{DHP_03}, \cite{Gu_02}, \cite{HHH_06}, \cite{MPS_00}, 
\cite{PS_05}, \cite{DK_11_1,DK_11_2,DK_11_3,DK_15},
and the references therein.

In what concerns the mixed norm results
with time dependent coefficients, for quite a long time
the power of summability with respect to the time variable
was assumed to be greater than that with respect
 to the space variables. Of course, in divergence form
equations such a restriction was not necessary due to
duality arguments. A remarkable step forward 
in this annoying problem
occurred in 2017
when the authors of \cite{GV_17_1} 
noticed that one can use the Muckenhoupt
$A_{p}$-weights and the Rubio de Francia
extrapolation theorem to get rid of the restriction
on the powers of summability even in the case of equations
with coefficients measurable in time variable.

Actually, the Muckenhoupt
weights were already used before for similar purposes.
It seems that \cite{HHH_03} was one of the first  (the authors of \cite{HHH_03}
cite the dissertation \cite{Fr_01} which appeared earlier
with similar approach)
  articles using the Muckenhoupt weights and the
Rubio de Francia extrapolation theorem
in the theory of evolution or elliptic equations
in spaces with mixed norms.
The authors deal with  higher-order parabolic equations with
{\em time independent continuous
coefficients\/} and use the approach based on the so-called
$R$-boundedness. 
In \cite{HH_03} the solvability is established
of the same type of
equations, with  time independent coefficients 
belonging to VMO, in spaces with $A_{p}$-weight as 
the underlying measures. Again the authors use
the Muckenhoupt weights and the
Rubio de Francia extrapolation theorem
to check the $R$-boundedness of certain families of operators.

In \cite{GV_17_1}  the authors consider 
parabolic higher-order 
equations with the coefficients
measurable in time and uniformly {\em continuous\/}
in the space variables.
They derive $L_{p}(L_{q})$-estimates for arbitrary
$p,q\in (1,\infty)$.
The $A_{p}$-weights and 
Rubio de Francia extrapolation theorem are used
again to check the  
$R$-boundedness of certain operators but also to 
use Mikhlin's multiplicator theorem and extrapolate
with respect to $t$, which was a natural
(for people who knew the extrapolation theorem) but 
absolutely new (for those, like the author,
 unfamiliar with weights
and extrapolation)
idea.
 The starting point in 
\cite{GV_17_1} is the case $q=p$ where
Mikhlin's theorem is used for equations
with constant coefficients.
The reader will find in \cite{ST_18}
a treatment of parabolic equations 
with coefficients depending only on time in mixed
norm spaces based on the Calder\'on-Zygmund theorem
rather than the Mikhlin theorem.
 
The results  of \cite{GV_17_1}
are generalized in \cite{GV_17_2}  
for parabolic systems satisfying the Legendre-Hadamard condition.

In the author's  opinion even more important step forward
was done in \cite{DK_18}, where, as in many cases, the general
ideology based on the Fefferman-Stein theorem (\cite{FS_72})
 is more or less standard
albeit its implementations vary. The
starting point
consists of   pointwise
estimates of the sharp functions of the derivatives
of solutions, which together with an appropriate
version of the Fefferman-Stein theorem (proved
in \cite{DK_18})
allows one to avoid using singular integrals,   the 
Mikhlin  or Calder\'on-Zygmund   theorems, $R$-boundedness, $H^{\infty}$-calculus,
and some other notions and tools from functional analysis.
In \cite{GV_17_1}, after a semigroup corresponding
to the operators with constant coefficient
is constructed in $L_{q}$-spaces with weights,
the authors use their special technique based on the fact that
$(-\Delta)^{m/2}$ has a bounded $H^{\infty}$-calculus
of angle $<\pi/2$ and $L_{q}(\bR^{d})$ has finite cotype
and also the 
$R$-boundedness of special kind of operators,
proved earlier in \cite{GLV_18}, to obtain their
main result.
In the proof of the $R$-boundedness they used the extrapolation
theorem.

The equations in \cite{DK_18} have the coefficients measurable
in $t$ and almost VMO in $x$, much more general than
in \cite{GV_17_1}. 
Furthermore, the authors of \cite{DK_18}
considered divergence and non-divergence equations whose coefficients are measurable
in one spatial variable.

The approach developed in \cite{DK_18} turned
out to be applicable even to fully nonlinear equations.
See \cite{DKr_18},
where the first to date weighted and
mixed-norm Sobolev estimates are presented
 for fully nonlinear elliptic
and parabolic equations in the whole space under a relaxed
convexity condition with almost VMO dependence on space-time
variables. The corresponding interior and boundary estimates
are also obtained.

The rest of the paper is organized as follows.
In the next section we present with almost all details
the proof of mixed-norm estimates for
the Laplacian. This proof is based
on the results in Section \ref{section 06.5.1.1}
about partitions and stopping times,
Section \ref{section 11.20.1}  about Muckenhoupt weights,
Section \ref{section 11.21.1} containing the proof
of Hardy-Littlewood maximal function theorem in $L_{p}$-spaces
with weights, Section \ref{section 11.21.10}  devoted
to the proof of the Rubio de Francia extrapolation
theorem, and Section \ref{section 11.29.1} about a generalized 
Fefferman-Stein theorem. In the  concluding
Section \ref{section 12.20.1} we state without proof
a result from recent paper \cite{DKr_18}
about mixed norm estimates for fully nonlinear equations.
 
\mysection{Illustration}
                                 \label{section 11.29.3}

Let 
$\bR^{d}$ be a $d$-dimensional Euclidean 
space with $d\geq2$ of points $x=(x_{1},...,x_{d})$.  
We consider Laplace's equation
\begin{equation}
                                          \label{11.29.5}
 \Delta u=f.
\end{equation}
 
Squaring   both sides of \eqref{11.29.5}
and integrating over $B_{r}:=\{x:|x|<r\}$ one easily shows
  (see, for instance,
Exercise 1.1.5 in \cite{Kr_08}) that, if $u$ is smooth
in $B_{r}$ and $u=0$ on $\partial B_{r}$,
then
\begin{equation}
                                       \label{11.30.1}
\sum_{i,j=1}^{d}\|D_{ij}u\|_{L_{2}(B_{r})}^{2}
\leq \|f\|_{L_{2}(B_{r})}^{2},
\end{equation}
where $D_{ij}=D_{i}D_{j}$, $D_{i}=\partial/\partial x_{i}$.
Then one proves (see, for instance, Exercise
1.3.23 in \cite{Kr_08}) that for any $f\in  L_{2}( B_{r} )$
there is a unique $u\in W^{2}_{2}( B_{r} )$ vanishing
on $\partial B_{r}$
such that $\Delta u=f$ in $B_{r}$. Moreover, it holds that
\begin{equation}
                                       \label{11.30.2} 
r^{-4}\|u\|^{2}_{ L_{2}(B_{r})}+r^{-2}
\sum_{i}\|D_{i}u \|^{2}_{ L_{2}(B_{r})}
+\sum_{i,j}\|D_{ij}u\|^{2}_
{ L_{2}(B_{r})}\leq 
5\|f\|^{2}_{ L_{2}(B_{r})},
\end{equation}
where $5$ is certainly not the best constant.

By using this fact and considering $u\zeta$, where
$\zeta$ is a cut-off function one proves (see
Exercise 2.4.6 of \cite{Kr_08} in which take $u$ independent
of $t$) that, if 
$u\in W^{2}_{2 }(B_{R})$ is harmonic in $B_{R}$,
then for any $r\in(0,R)$
\begin{equation}
                                       \label{12.2.6} 
\|D^{2}u\|_{L_{2}(B_{r})}\leq N(d)(R-r)^{-2}\|u\|_{L_{2}(B_{R})},
\end{equation}
where and below by $D^{m}u$ we mean the collection
of all $m$-th order partial derivatives of $u$ and by
it $L_{p}$-norm we mean the $L_{p}$-norm of
$|D^{m}u|$, the latter being the Euclidean norm
of $D^{m}u$ as a vector in an appropriate Euclidean
space. 
One can iterate this estimate by applying it
to the derivatives of $u$, which are harmonic
in $B_{R}$ along with $u$. Then one sees that,
for any $m\geq 1$,
$$
\| D^{2m}u\|_{L_{2}(B_{r})}\leq
 N(d,r,R,m)\|D^{2}u\|_{L_{2}(B_{R})}.
$$
For $m$ large enough, by Sobolev embedding theorems
it follows that
\begin{equation}
                                       \label{11.30.3}
\max_{B_{1}}|D^{3}u|\leq N(d)\|D^{2}u\|_{L_{2}(B_{2})}
\end{equation}
  if 
$u\in W^{2}_{2 }(B_{2})$ is harmonic in $B_{2}$.

By using the notation  
$$
\dashint_{A}f\,dx=\frac{1}{|A|}\int_{A}f\,dx,
$$
where $|A|$ is the volume of $A$,
we infer that for the functions $u$ as in \eqref{11.30.3}
and $r\in(0,1]$
we have
$$
\dashint_{B_{r}}\dashint_{B_{r}}|D^{2}u(x)
-D^{2}u(y)|^{2}\,dxdy\leq r^{2}N(d)\int_{B_{2}}
|D^{2}u|^{2}\,dx.
$$
By using scaling we obtain the following important
estimate (we have just repeated in the simplest situation
the proof of Theorem 4.2.6 of \cite{Kr_08}).

\begin{lemma}
                                    \label{lemma 11.30.2}
Let $\nu\geq 2$ and $r\in(0,\infty)$ be some
constants and let $u\in W^{2}_{2}(B_{\nu r})$
be harmonic in $B_{\nu r}$. Then with a constant $N=N(d)$
we have
\begin{equation}
                                       \label{11.30.4}
\dashint_{B_{r}}\dashint_{B_{r}}|D^{2}u(x)
-D^{2}u(y)|^{2}\,dxdy\leq N\nu^{-2} \dashint_{B_{\nu r}}
|D^{2}u|^{2}\,dx.
\end{equation}
\end{lemma}

Now we have an analog of Theorem 4.3.1 of \cite{Kr_08}.
\begin{theorem}
                                   \label{theorem 11.30.1} 
Let $\nu\geq 2$ and $r\in(0,\infty)$ be some
constants and let $u\in W^{2}_{2}(B_{\nu r})$. 
Then with a constant $N=N(d)$
we have
$$
\dashint_{B_{r}}\dashint_{B_{r}}|D^{2}u(x)
-D^{2}u(y)|^{2}\,dxdy
$$
\begin{equation}
                                       \label{11.30.40}
\leq N\nu^{d }\dashint_{B_{\nu r}}
|\Delta u|^{2}\,dx +
N\nu^{-2} \dashint_{B_{\nu r}}
|D^{2}u|^{2}\,dx.
\end{equation}
\end{theorem}

Proof. Let $v\in W^{2}_{2}(B_{\nu r})$
be the solution of $\Delta v=\Delta u$,
that equals zero on $\partial B_{\nu r}$.
By taking into account \eqref{11.30.2}, we get
$$
\dashint_{B_{r}}\dashint_{B_{r}}|D^{2}v(x)
-D^{2}v(y)|^{2}\,dxdy\leq 
4\dashint_{B_{r}}|D^{2}v |^{2}\,dx
$$
\begin{equation}
                                       \label{11.30.50}
\leq 4\nu^{d}\dashint_{B_{\nu r}}|D^{2}v |^{2}\,dx
\leq N\nu^{d} \dashint_{B_{\nu r}}
|\Delta u|^{2}\,dx.
\end{equation}
 
 On the other  hand, $w:=u-v$ is harmonic and by Lemma
\ref{lemma 11.30.2}
$$
\dashint_{B_{r}}\dashint_{B_{r}}|D^{2}w(x)
-D^{2}w(y)|^{2}\,dxdy\leq N\nu^{-2} \dashint_{B_{\nu r}}
|D^{2}w|^{2}\,dx
$$
$$
\leq N\nu^{-2} \dashint_{B_{\nu r}}
|D^{2}u|^{2}\,dx+N\nu^{-2} \dashint_{B_{\nu r}}
|D^{2}v|^{2}\,dx,
$$
which combined with \eqref{11.30.50} leads 
to \eqref{11.30.40} and proves the theorem.

Then let $\bB$ be the collection of open balls in
$\bR^{d}$ and introduce the notation  
$$
f^{\#}(x)=\sup_{ B\in\bB,B\ni x }
\dashint_{B_{r}}\dashint_{B_{r}}|f(x')
-f(x'')| \,dx'dx'',
$$
$$
\bM f(x)=\sup_{  B\in\bB,B\ni x } 
\dashint_{B_{r}} |f(y)
 | \,dy .
$$
The function $f^{\#}$ is called the sharp function of
$f$ and $\bM f$ is its maximal function.

In light of H\"older's inequality and
 the possibility of changing the origin
Theorem \ref{theorem 11.30.1} implies the following
estimate, the direct analogs of which play  a crucial role
in the theory of linear and fully nonlinear equations
with VMO main coefficients.

\begin{theorem}
                                   \label{theorem 12.1.1}
There is a   constant $N=N(d)$ such that for any 
  $\nu\geq 2$  
  and  $u\in W^{2}_{2,\loc}(\bR^{d})$
we have
\begin{equation}
                                       \label{11.30.04}
(D^{2}u)^{\#}\leq N\nu^{d /2}(\bM(
|\Delta u|^{2}))^{1/2}  +
N\nu^{-1} (\bM(
|D^{2}u|^{2}))^{1/2}.
\end{equation}
\end{theorem}
The approach based on pointwise estimates involving sharp and
maximal functions 
was first suggested in \cite{Kr_07}.

We can raise both sides
of \eqref{11.30.04} to   the power $p>2$ and then integrate through 
 by using any measure $w(dx)$,
but, if this measure has a density $w\in A_{p/2}$
(see the definition of $A_{p}$ in Section \ref{section 11.20.1}),
 then,
on the one hand,
the Hardy-Littlewood theorem turns out to hold for such a measure
(see Theorem \ref{theorem 11.21.3}) and one can estimate
the integrals of $(\bM(
|\Delta u|^{2}))^{p/2}$ and $(\bM(
|D^{2}u|^{2}))^{p/2}$ just with the integrals of 
$|\Delta u|^{p}$ and $|D^{2}u|^{p} $ against $w(dx)$,
respectively. On the other hand the Feffeman-Stein theorem
turns out to be true (see Theorem \ref{theorem 12.19.1}), 
which allows one
to estimate the integral of $|D^{2}u|^{p}$ through
the integral of  the $p$th power of its sharp function.
Finally, having $\nu^{-1}$ on the right allows us to absorb
what is coming from the last term in \eqref{11.30.04}
into the integral of $|D^{2}u|^{p}$ (provided that it is finite)
and leads to the following in which $[w]_{A_{p/2}}$
is the $A_{p/2}$-constant of the $A_{p/2}$-weight $w$,
which affects  both constants in
the Hardy-Littlewood   and the Feffeman-Stein theorems.

\begin{theorem}
                                   \label{theorem 12.1.2}
Let $p>2$ and $w$ be an $A_{p/2}$-weight on $\bR^{d}$. Then
there  is a   constant $N=N(d,p,[w]_{A_{p/2}})$ such that for any 
  $u\in W^{2}_{p}(\bR^{d})$
we have
\begin{equation}
                                       \label{11.30.5}
\int_{\bR^{d}}|D^{2}u|^{p}w\,dx\leq N
\int_{\bR^{d}}|\Delta u|^{p}w\,dx.
\end{equation}
\end{theorem}

Of course, raising \eqref{11.30.04} to any power $p>2$
and applying the usual versions of 
the Hardy-Littlewood   and the Feffeman-Stein theorems
rather than their weighted counterparts leads to the classical
estimate 
\begin{equation}
                                       \label{11.30.6}
\|D^{2}u\|_{L_{p}}\leq N(p,d)\|\Delta u\|_{L_{p}},
\quad \forall u\in W^{2}_{p}(\bR^{d}),
\end{equation}
first obtained by using the Calder\'on-Zygmund theorem.
This is what was usually done in parts of the theory
related to VMO conditions and using the Fefferman-Stein theorem.

What is amazing is that \eqref{11.30.5}
 implies not only \eqref{11.30.6}
for $p>2$, but also the following
result
 by   just a mere reference to a simple corollary
(see Theorem \ref{theorem 12.2.1})
of the Rubio de Francia extrapolation theorem
(Theorem \ref{theorem 11.21.6}).
Note that Theorem \ref{theorem 12.1.3},
  is a very particular case of the results
in \cite{DK_18}.

To state this result,
for $i\in\{1,...,d\}$ and
 $p_{i}\in(1,\infty)$,   introduce
$$
\| f
\|_{L_{p_{1},\ldots,p_{d}} }^{p_{d}}
$$
\begin{equation}
                                       \label{4.9.6}
:=\int_{\bR}\Big(\cdots
\Big(\int_{\bR}
\Big(\int_{\bR } |f|^{p_{1}} \,d  x_{1}
\Big)^{p_{2} /p_1 } \,d  x_{2}\Big)^{p_{3}/p_{2}}
\cdots\Big)^{p_{d} /p_{d-1} } \,d  x_{d}.
\end{equation}

\begin{theorem}
                                   \label{theorem 12.1.3}
Let $p_{1},...,p_{d}\in(2,\infty)$. Then
there  is a   constant $N=N(d,p_{1},...,p_{d})$ such that for any 
  $u\in W^{2}_{2,\loc }(\bR^{d})$
we have
\begin{equation}        \label{11.30.7}
\| D^{2}u
\|_{L_{p_{1},\ldots,p_{d}} } \leq N
\| \Delta u
\|_{L_{p_{1},\ldots,p_{d}} } ,
\end{equation}
provided that $\| D^{2}u
\|_{L_{p_{1},\ldots,p_{d}} } +\| D u
\|_{L_{p_{1},\ldots,p_{d}} }+\|  u
\|_{L_{p_{1},\ldots,p_{d}} } $ is finite. 
\end{theorem}
\begin{remark}
                                     \label{remark 12.1.3}
This ``provided that ... is finite''
appears because we obtained \eqref{11.30.5} for
$u\in W^{2}_{2}(\bR^{d})$ rather than 
$u\in W^{2}_{2,\loc }(\bR^{d})$ and the derivation 
of \eqref{11.30.7} first proceeds for $u\zeta$, where
$\zeta$ is a cut-off function and then sending
$\zeta\to1$. Similarly to the case of \eqref{11.30.6},
\eqref{11.30.7} is false if we drop the above provision.
The example $u=x_{1}^{2}-x_{2}^{2}$ shows this.
\end{remark}

To derive some direct
consequences of Theorem \ref{theorem 12.1.3} consider
equation \eqref{11.29.5} in $\bR^{d}_{+}=\{x=(x_{1},x'):
x_{1}\in[0,\infty),x'\in\bR^{d-1}\}$ with either
zero Dirichlet or zero Neumann condition. By using
either odd or even continuation across the plane
$\{x_{1}=0\}$ one obtains the following estimates
for $p_{1},p_{2}>2$ in both cases:
\begin{equation}
                                       \label{12.2.2}
\int_{\bR^{d-1}}\Big(\int_{0}^{\infty}|D^{2}u|^{p_{1}}
\,dx_{1}\Big)^{p_{2}/p_{1}}\,dx'
\leq N\int_{\bR^{d-1}}\Big(\int_{0}^{\infty}|\Delta u|^{p_{1}}
\,dx_{1}\Big)^{p_{2}/p_{1}}\,dx',
\end{equation}
\begin{equation}
                                       \label{12.2.3}
\int_{0}^{\infty}\Big(\int_{\bR^{d-1}}|D^{2}u|^{p_{2}}
\,dx'\Big)^{p_{1}/p_{2}}\,dx_{1}
\leq N\int_{0}^{\infty }\Big(\int_{\bR^{d-1}}|\Delta u|^{p_{2}}
\,dx'\Big)^{p_{1}/p_{2}}\,dx_{1}.
\end{equation}
Of course, these estimates are only true
provided that conditions similar to the one in
Theorem \ref{theorem 12.1.3} are satisfied.

Above one can also replace $dx_{1}$ with $w(x_{1})\,dx_{1}$,
where $w$ is any $A_{p_{1}/2}$-weight
(see Theorem \ref{theorem 12.2.1}) and then the constants
$N$ depend also on $[w]_{A_{p_{1}/2}}$. By the way,
we know that
in $(0,\infty)$ the functions $x ^{q}$ are $A_{p}$-weights
iff $q\in(-1,p-1)$.

The outlined way of proving  
\eqref{11.30.6}
is similar to what is done in the
part of the Sobolev-space theory
of linear and fully nonlinear elliptic and parabolic
equations which rely on the Fefferman-Stein theorem.
There is another part where they use an approach based
on the theory of viscosity solutions and which is
out of the scope of the present review. Still
it would be interesting to know if one can obtain,
say  Theorem \ref{theorem 5.1.1}, by using
the methods of the theory of viscosity solutions.

One may  be not satisfied with the restriction $p_{i}>2$
in \eqref{11.30.7}, \eqref{12.2.2}, and \eqref{12.2.3}.
There are two known ways around it:

(a) Take $p _{0}\in(1,2)$ and first prove (or use) the fact that
for any $f\in L_{p}(B_{r})$
equation \eqref{11.29.5} in $B_{r}$  with zero
boundary condition on $\partial B_{r}$ has a unique
solution $u\in W^{2}_{p}(B_{r})$ and a natural analog
of \eqref{11.30.2} holds with $L_{p}(B_{r})$
in place of $L_{2}(B_{r})$. Then by just repeating
what is below \eqref{11.30.2} we come to the 
conclusion that \eqref{11.30.7}, \eqref{12.2.2},
 and \eqref{12.2.3}
hold for any $p_{i}\in (p_{0},\infty)$, and since $p_{0}>1$
is arbitrary, they hold for any $p_{1},p_{2}\in (1,\infty)$.

(b) Take $p _{0}\in(1,2)$ and use \eqref{11.30.6} 
(which, as we pointed out, follows from
Theorem \ref{theorem 12.1.1}  for $p>2$ and then for $p\in(1,2)$
by duality) to derive \eqref{12.2.6} 
for $u\in W^{2}_{p_{0}}(B_{R})$
with $L_{p}$ in place of $L_{2}$ as in
Exercise 2.4.6 of \cite{Kr_08} in which take $u$ independent
of $t$. This will lead to Lemma \ref{lemma 11.30.2}
for the same range of $\nu$ with power of summability $2$
replaced by $p_{0}$ and then lead to Theorem
\ref{theorem 11.30.1} with the same replacement
of $2$ by $p$ but for $\nu\geq 4$. This 
new restriction comes from the fact that this time we do
not want to use the solvability of \eqref{11.29.5} in balls
and instead take a $\zeta\in C^{\infty}_{0}(B_{\nu r})$
such that $\zeta=1$ on $B_{\nu r/2}$, solve the equation
$\Delta v=\zeta\Delta u$ in $\bR^{d}$ and define
$w=u-v$ which turns out to be harmonic in $B_{\nu r/2}$.
Then an obvious modification of the short proof
of Theorem \ref{theorem 11.30.1} proves it with $p_{0}$
in place of $2$. After that one proceeds as in (a). 
This way is used even in the first step
for equations with VMO coefficients when
the solvability of equations in  question
in balls or cylinders is either unknown or
hard to obtain, in particular, because of  the 
necessity to deal with 
boundary conditions (higher-order case).

\mysection{Partitions and stopping times}
                      \label{section 06.5.1.1}

Fix some integers $k_{1},...,k_{d}\geq1$
and call any 
$$
C_{l,x}=x+[0,l ^{ k_{1}} )\times
...\times [0,l ^{ k_{d} }),\quad x\in\bR^{d},l>0
$$
 a ``cube'' or half-closed ``cube''
of size $l$.  
Observe that if $k_{1}=...=k_{d}=1$, then we deal just 
with usual cubes whose edges are parallel to the
coordinate axes and which are common
in the theory of elliptic equations. If $k_{1}=2$ and 
$k_{2}=...=k_{d}=1$,
we are dealing with parabolic ``cubes'' often
used in the theory of second-order parabolic equations.

Next, we introduce a subset $\Omega$ of $\bR^{d}$
which is, actually a product of some $\bR$
and some shifted $[0,\infty)$. Namely, we assume that
we are given $\Omega\subset\bR^{d}$ such that,
for each $l>0$, it is the disjoint union of some
``cubes'' of size $l$ which belong to $\Omega$.
The simplest example is $\Omega=\bR^{d}$.
The reader is advised to always have in mind this basic 
example.

For $n\in\bZ=\{0,\pm1,...\}$
introduce the following families of ``dyadic cubes''
by  
$$
\bC_{n}=\{C_{n}(i_{1},...,i_{d}):
i_{1},...,i_{d}\in\bZ,C_{n}(i_{1},...,i_{d})
\subset\Omega\},
$$
 where
$$
C_{n}(i_{1},...,i_{d})=[i_{1}2^{-k_{1}n},(i_{1}+1)2^{-k_{1}n})\times...\times
[i_{d}2^{-k_{d}n},(i_{d}+1)2^{-k_{d}n}).
$$
Also set $\bC=\bigcup_{n\in \bZ}\bC_{n}$.

Define $\cF_{n}$ as the collection of subsets of $\Omega$
consisting
of an empty set and of the unions
of some elements of $\bC_{n}$.
Obviously $\cF_{n}\subset\cF_{m}$
for $n\leq m$. 

If $\tau=\tau(x)$ is a function on $\Omega$ with values in
$\{\infty,0,\pm1,\pm2,...\}$, we call $\tau$ {\em a stopping time\/}
(relative to $\{\cF_{n}\}$)
if, for each $n=0,\pm1,\pm2,...$,  
$$
\{x:\tau(x)=n\}\in\cF_{n}
$$ 
(that is $\{x:\tau(x)=n\}$ is
 either empty or else is the union of
some elements of $\bC_{n}$).
The simplest example of a stopping time is given by
$\tau(x)\equiv0$. 

If $\tau$ is a stopping time we denote by $\cF_{\tau}$
the collection of 
  Borel $A$ such that, for any $n\in\bZ$
we have
$$
A\cap\{\tau=n\}\in \cF_{n}.
$$
Observe that
if we are given two stopping times $\tau$ and $\sigma$
and $\sigma\leq \tau$, then $\cF_{\sigma}
\subset\cF_{\tau}$ since 
$$
 A\cap\{\tau=n\}=\bigcup_{k=-\infty}^{n }
 A\cap\{\sigma=k\}\cap\{\tau=n\}
$$
and $A\cap\{\sigma=k\}\in\cF_{k}\subset \cF_{n}$.
Obviously the intersection of two sets in $\cF_{\tau}$
belongs to $\cF_{\tau}$. An easy and useful
fact is that $\Omega,\{\tau<\infty\}\in \cF_{\tau}$.
Also a useful fact to remember is that if $A\in\cF_{\tau}$,
then $A\cap\{\tau<\infty\}$ is a disjoint union
of $A\cap\{\tau=n\}\in \cF_{n}$ and
each $A\cap\{\tau=n\}$ is either empty or is a disjoint union of
some $C\in\bC_{n}$
such that $C=C\cap\{\tau=n\}$, the latter showing that
$C\in \cF_{\tau}$.
 
We assume that we are given a measure $\mu$
on Borel subsets of $\Omega$ such that, for any $x\in\Omega$,
 \begin{equation}
                                                  \label{11.28.1}
\lim_{l\to\infty} \mu(C_{l,x} )=\infty.
\end{equation}

Whenever it makes sense we
  use the notation 
$$
  f_{A}=\dashint_{A}f\,\mu(dx):=\frac{1}{\mu(A)}\int_{A}f(x)\,\mu(dx)  
\quad \bigg(\frac{0}{0}:=0\bigg)
$$
 for  the average value of $f$ over $A$. 
 
Next, for each $x\in \Omega$ and $n\in \bZ$ 
there exists (a unique)
$C\in\bC_{n}$ such that $x\in C$. We denote 
this $C$ by $C_{n}(x)$.
The reader can check himself (or consult,
for instance, \cite{Kr_08} for this and a few more simple 
properties of introduced objects) that
for any $\lambda>0$ and $f\geq 0$, such that
$ f _{C_{n}(x)}\to0$ as $n\to-\infty$
for any $x$, 
\begin{equation}
                                            \label{11.20.1}
\tau_{\lambda}(x)=\inf\{n: f_{C_{n}(x)}>\lambda\}
\quad (\inf\emptyset:=\infty)
\end{equation}
is a stopping time.

 We also use the notation
$$
f_{|n}(x)= f_{C_{n}(x)}=\dashint_{C_{n}(x)}f(y)\,\mu(dy).
$$
If we are also given a stopping time $\tau$, we let
$$
f_{|\tau}(x) =f_{|\tau(x)}(x)
$$
 for those $x$ for which
$\tau(x)<\infty$  and $f_{|\tau}(x)=f(x)$ otherwise.

We suppose that $\mu$ satisfies the ``doubling condition'':
  for any $n$, $C\in \bC_{n}$, and $C'\in\bC_{n-1}$ such that 
$C\subset C'$ we have
\begin{equation}
                                            \label{11.21.2}
\mu(C')\leq N_{0}\mu(C),
\end{equation}
 where $N_{0}$
is a constant independent of $n,C, C'$.
One of consequences of this condition is that
for $f\geq0$
on the set where $\tau_{\lambda}(x)=n$
we have $f_{|n-1}(x)\leq\lambda$ and
\begin{equation}
                                     \label{11.23.1}
f_{\tau}(x)=\frac{1}{\mu(C_{n}(x))}
\int_{C_{n}(x)}f\,\mu(dy)
\leq \frac{N_{0}}{\mu(C_{n-1}(x))}
\int_{C_{n-1}(x)}f\,\mu(dy)\leq N_{0}\lambda.
\end{equation}
 
Another consequence of \eqref{11.28.1} and \eqref{11.21.2}
is that $\mu(C_{l,x})>0$ for any $x\in \Omega$ and $l>0$.
 
In the following lemma by $I_{A,\tau<\infty}$
we mean the indicator function of the set
$\{x\in A:\tau(x)<\infty\}$.
\begin{lemma}
                                 \label{lemma 06.5.30.1}

(i) Let   $f$ be Borel on $\Omega$, $f\geq0$,
  let $\tau$ be a stopping time, and let $A\in \cF_{\tau}$.  
  Then
\begin{equation}
                                          \label{06.5.1.1}
\int_{\Omega} f_{|\tau}(x) I_{A,\tau<\infty}\,\mu(dx)=
\int_{\Omega} f(x)  I_{A,\tau<\infty}\,\mu(dx) .
\end{equation}

(ii) Let   $f$ be Borel  on $\Omega$, $f\geq0$,
and let $\lambda>0$ be a constant. 
Assume that $f_{|n}(x)\to0$ as $n\to-\infty$
at any $x$.
Then for $\tau=\tau_{\lambda}$ defined in \eqref{11.20.1} 
we have
\begin{equation}
                                          \label{11.18.6}
\lambda I_{\tau<\infty} <f_{|\tau}(x)
I_{\tau<\infty}\leq N_{0}\lambda,
\end{equation}
and for any $A\in\cF_{\tau}$
$$
N_{0}^{-1}\lambda^{-1}\int_{\Omega}f(x)
I_{A,\tau<\infty}\,\mu(dx)\leq 
\mu(\{x\in A:\tau(x)<\infty\})
$$
\begin{equation}
                                          \label{06.5.1.3}
\leq 
 \lambda^{-1}\int_{\Omega}f(x)
I_{A,\tau<\infty}\,\mu(dx).
\end{equation}

\end{lemma}

 Proof. (i) Owing to the additivity of
the integral, it suffices to prove
\eqref{06.5.1.1} with $\tau=n$ in place
of $\tau<\infty$ and, since the set
$\{x\in A:\tau(x)=n\}$ is the disjoint
union of some $C\in \bC_{n}$, it only remains to observe that
for such $C$ we have $C\cap\{\tau=n\}=C$ and
$$
\int_{\Omega} f_{|\tau}(x) I_{C,\tau=n}\,\mu(dx)
 =\int_{\Omega} f_{C} I_{C}\,\mu(dx)
=\int_{C}f\,\mu(dx)=\int_{\Omega}fI_{C,\tau=n}\,\mu(dx).
$$

(ii) Relations \eqref{11.18.6} follow
from the definition of $\tau$ and \eqref{11.23.1}.
The first inequality in \eqref{06.5.1.3}
follows from \eqref{06.5.1.1} and \eqref{11.18.6}.
The second one follows from Chebyshev's
inequality and \eqref{06.5.1.1} because
$$
\mu(\{x\in A:\tau(x)<\infty\})=
\mu(\{x\in \Omega:f_{\tau}I_{A,\tau <\infty}>\lambda\}).
$$
The lemma is proved.

Define
{\em the maximal ``dyadic'' function\/}
of $f$ by
\begin{equation}
                           \label{12.28.1}
\cM f(x)=\sup_{n<\infty}|f|_{|n}(x),
\end{equation}
 so that
$\cM f=\cM |f|$.
 Observe that, if $f\geq0$, $\{x:\cM f(x)>\lambda\}=\{x:\tau_{\lambda}<\infty\}$,
where $\tau_{\lambda}$ is taken from
\eqref{11.20.1}. Therefore, Lemma \ref{lemma 06.5.30.1}  (ii) implies the following.

\begin{corollary}
                    \label{corollary 12.25.1}
Under the conditions of Lemma \ref{lemma 06.5.30.1}  (ii)
$$
N_{0}^{-1}\lambda^{-1}\int_{\Omega}f(x)
I_{\cM f>\lambda}\,\mu(dx)\leq 
\mu(\{x\in A:\cM f(x)>\lambda\})
$$
\begin{equation}
                       \label{12.25.2}
\leq 
 \lambda^{-1}\int_{\Omega}f(x)
I_{\cM f>\lambda}\,\mu(dx).
\end{equation}
\end{corollary}

The following standard consequence of
the second inequality in \eqref{12.25.2}
is left to the reader as an exercise
(see Exercise 3.2.7 and the hint to it in \cite{Kr_08}).
\begin{corollary}
                    \label{corollary 12.25.2}
Let   $fI_{C}\in \cL_{1}(\mu)$
for any $C\in\bC$. Then
(Lebesgue differentiating theorem)
$f_{|n}\to f$ ($\mu$-a.e.)
and, in particular, 
$|f|\leq \cM f$ ($\mu$-a.e.). 
\end{corollary}
It is worth noting that both inequalities in
\eqref{06.5.1.3} are crucial
in the proof of the reverse H\"older's inequality
for Muckenhoupt's weights, and the right
inequality is used in a crucial way in proving
the Fefferman-Stein theorem in Sect on \ref{section 11.29.1}.

The following remark will not be used in the future.
It was hard not to make it. 
\begin{remark}
                                 \label{remark 12.20.1}
The first inequality in \eqref{06.5.1.3} is instrumental
not only in proving
  the reverse H\"older
inequality for $A_{p}$-weights but its version
\eqref{12.25.2}
also
is crucial  in the proof of the first part
of a remarkable
Zygmund-Stein result that, for any 
Borel $f\geq0$ and $\lambda_{0}>0$,
\begin{equation}
                                      \label{12.20.3}
\int_{\Omega}\cM fI_{\cM f>\lambda_{0}}\,\mu(dx)<\infty
\Longrightarrow 
\int_{\Omega}  fI_{  f>\lambda_{0}}\log(f/\lambda_{0})\,\mu(dx)<\infty,  
\end{equation}
\begin{equation}
                                      \label{12.20.4}
\int_{\Omega}  fI_{  f>\lambda_{0}}\log(f/\lambda_{0})\,\mu(dx)<\infty
\Longrightarrow 
\int_{\Omega}\cM fI_{\cM f>2\lambda_{0}}\,\mu(dx)<\infty.
\end{equation}

Here \eqref{12.20.3} is obtained just by integrating
with respect $\lambda\in(\lambda_{0},\infty)$   the
first inequality in \eqref{12.25.2}, where on the left
$\cM f$ is replace with a smaller quantity $f$.
To prove \eqref{12.20.4}   use
the second inequality in \eqref{12.25.2}, which
after observing that $\cM (fI_{ f \leq \lambda})\leq\lambda$
implies that
$$
  \mu(
\{x:\cM f(x)>2\lambda\})\leq
 \mu(
\{x:\cM (fI_{ f >\lambda})(x)> \lambda\})
$$
$$
\leq \lambda^{-1}\int_{\Omega}
f(x)I_{ f >\lambda} \,\mu(dx).
$$
Then again integrate
with respect $\lambda\in(\lambda_{0},\infty)$ the inequality between the extreme terms.
After that one will only need to estimate $2
\lambda_{0}\mu(
\{x:\cM f(x)>2\lambda_{0}\})$  and observe that
this quantity is less than
$$
2\lambda_{0}
\mu(
\{x:\cM fI_{f>3\lambda_{0} /2}(x)>3\lambda_{0} /2\})
\leq(4/3)\int_{\Omega}
fI_{f>3\lambda_{0} /2} \,\mu(dx),
$$
where the last term is finite due to  the condition in \eqref{12.20.4}.

\end{remark}

\mysection{Muckenhoupt's weights}
                                 \label{section 11.20.1}  
We use the setting of Section
\ref{section 06.5.1.1}, but the doubling condition
here is stronger than \eqref{11.21.2}.
We need \eqref{11.21.2} to hold for 
a collection of subsets of $\Omega$ which contains
 all translates of
``contracted'' interiors of
 $C$ which keep them in $\Omega$.
For $x\in \Omega$ and $l>0$ denote
$$
D_{l,x}=x+(0,l ^{ k_{1}} )\times
...\times (0,l ^{ k_{d} }).
$$
Assume that we have a family $\bQ$ of open subsets
of $\Omega$ which contains all $D_{l,x}$ and is such that
for any $Q\in\bQ$ there exists a $D_{l,x}$ satisfying
\begin{equation}
                                   \label{11.24.1}
Q\subset D_{l,x},\quad \mu(D_{2l,x})\leq N_{0}
\mu(Q).
\end{equation}

We suppose that 
\begin{equation}
                                   \label{11.24.4} 
\mu(\partial \Omega)=0
\end{equation}
and then, as is easy to see, the 
new doubling condition 
\eqref{11.24.1} implies \eqref{11.21.2} and, moreover,
$\mu(C_{2l,x})\leq N_{0}\mu(C_{l,x})$
for any $x\in\Omega$ and $l>0$
(just in case, observe that the sets $D_{l,x}$
are open).

\begin{definition}[Muckenhoupt's weights]
                              \label{definition 11.17.1}
Let $w(x)$ be a  function
on $\Omega$ such that $0< w<\infty$ ($\mu$-a.e.)
and $w_{Q}<\infty$ for any $Q\in\bQ$.
We call it an $A_{1}$-weight 
(relative to $\bQ$ or relative to $(\bQ,\mu)$) if there is a constant
$N$ such that
\begin{equation}
                                    \label{11.19.1}
w_{Q}\leq Nw(x)\quad\forall x\in Q, \forall Q\in\bQ.
\end{equation}
on $\Omega$. The least constant $N$ satisfying 
\eqref{11.19.1} is called the $A_{1}$-constant of $w$
  denoted by $[w]_{A_{1}}$. For $p\in(1,\infty)$
we call $w$ an $A_{p}$-weight if there is a constant
$N$ such that
\begin{equation}
                                    \label{11.17.2}
w_{Q}\Big(\big(w^{-1/(p-1)}\big)_{Q}\Big)^{p-1}\leq N\quad \forall Q\in\bQ.
\end{equation}
  The least constant $N$ satisfying 
\eqref{11.17.2} is called the $A_{p}$-constant of $w$
 denoted by $[w]_{A_{p}}$.

\end{definition}
 
It is known that, if $\Omega=\bR^{d}$
with Lebesgue measure, $|x|^{\alpha}$ is an $A_{p}$-weights 
($p\in[1,\infty)$)
 if and only if $-d<\alpha<(p-1)d$.

Set 
$$
 w(A)=\int_{A}w\,\mu(dx), 
$$
and denote $L_{p}(w)=L_{p}(\Omega,w(dx))$.
\begin{remark}
                                       \label{remark 11.25.1}
H\"older's inequality implies that, if $p\in(1,\infty)$,
   $A $ is Borel, and $\infty>w>0$ ($\mu$-a.e.) on $A$, then
\begin{equation}
                                    \label{11.25.2}
1=\dashint_{A}w^{1/p}w^{-1/p}\,\mu(dx)
\leq (w_{A})^{1/p}\Big(\big(w^{-1/(p-1)}\big)_{A}\Big)^{(p-1)/p}.
\end{equation}
Therefore, $[w]_{A_{p}}\geq1$ if $w\in A_{p}$.
 This  obviously holds
for $p=1$ as well.
 
Also note that one can replace what is raised to the
power $p-1$ in \eqref{11.17.2} with
$\big(w^{-1/(p-1)}\big)_{A}\mu(A)/\mu(Q)$
for any Borel $A\subset Q$. Then \eqref{11.25.2}
implies that 
\begin{equation}
                                    \label{11.17.4}
[w]_{A_{p}}\frac{w(A)}{w(Q)}
\geq\Big(\frac{\mu(A)}{\mu(Q)}\Big)^{p}.
\end{equation}

We obtained this result for $p>1$. It is also true
for $p=1$ since, for $w\in A_{1}$ we have $[w]_{A_{1}}
w\geq w_{Q}$
on $Q$ and hence
$$
\mu(A)=\int_{A}\frac{1}{w}w\,\mu(dx)
\leq [w]_{A_{1}}\frac{1}{w_{Q}}w(A),
$$
which is \eqref{11.17.4} for $p=1$.

An important consequence of \eqref{11.17.4}
is that for $p\in[1,\infty)$, $w\in A_{p}$, and
 any $\alpha\in(0,1)$ there exists $\beta\in(0,1)$
depending only on $\alpha$, $p$, and $[w]_{A_{p}}$, such that, 
for any   Borel $S\subset Q\in\bQ$ 
\begin{equation}
                                    \label{11.24.6}
\mu(S)\leq\alpha\mu(Q)\Longrightarrow
 w(S)\leq\beta w(Q).
\end{equation}
On gets \eqref{11.24.6} by inspecting the following
version of \eqref{11.17.4}
$$
[w]_{A_{p}}\Big(1-\frac{w(S)}{w(Q)}\Big)
\geq\Big(1-\frac{\mu(S)}{\mu(Q)}\Big)^{p}\geq (1-\alpha)^{p}.
$$
\end{remark}

\begin{remark}
                                      \label{remark 12.22.1}
In the future we several times say that
a constant entering an estimate
 depends only on ..., and $[w]_{A_{p}}$.
It is important to emphasize that, given a 
$K_{0}$, in these situations
the constant can be chosen  to
depend  only on ..., and $K_{0}$, provided
that $[w]_{A_{p}}\leq K_{0}$.

\end{remark}

As a simple exercise based on \eqref{11.19.1}
or H\"older's inequality one proves that
$A_{p}\subset A_{q}$ if $1\leq p\leq q$.
Here is an extension of this for $A_{1}$-weights
(when $v\equiv1$).

\begin{lemma}
                                 \label{lemma 11.19.1}
If $w ,v\in A_{1}$ and $p\in(1,\infty)$, then
$w v^{1-p}\in A_{p}$ and
$[w v^{1-p}]_{A_{p}}\leq[w] _{A_{1}}[v]^{p-1}_{A_{1}}$.
\end{lemma}

Proof. In light of \eqref{11.19.1} for any $Q\in\bQ$
$$
(w v ^{1-p})_{Q}\leq [v]_{A_{1}}^{p-1}(w (v_{Q})^{1-p})_{Q}
=[v]_{A_{1}}^{p-1}w_{Q}(v_{Q})^{1-p}.
$$
$$
(w^{-1/(p-1)}v)_{Q}\leq [w]_{A_{1}}^{ 1/(p-1)}
((w_{Q})^{ -1/(p-1)}v)_{Q}
=[w]_{A_{1}}^{ 1/(p-1)}(w_{Q})^{-1/(p-1)}v_{Q},
$$
and our assertions follow. The lemma is proved.

A deep result, which however we will not need,
is that, under very general conditions,
 any $A_{p}$-weight admits a representation
as $w v^{1-p}$ with $w,v\in A_{1}$
(P. Jones' theorem).

\begin{theorem}[Reverse H\"older's inequality]
                                 \label{theorem 11.17.3}
If $p\in[1,\infty)$ and $w \in A_{p}$, then
there exist constants
$N=N(p,N_{0},[w]_{A_{p}})<\infty$ and
 $\varepsilon=\varepsilon(p,N_{0},[w]_{A_{p}})>0$
such that
for any  
 $Q \in  \bQ $  we have
\begin{equation}
                                          \label{11.25.4}
  (w^{1+\varepsilon})_{Q}\leq N
(w_{Q})^{1+\varepsilon}.
\end{equation}
\end{theorem}

Proof. Since $A_{1}\subset A_{p}$
for any $p>1$, we may assume that $p>1$.
We also may assume that $ w(Q)>0$. Then we find
a $D_{l,x}$ such that \eqref{11.24.1} holds and observe
that, since $\mu(Q)\geq N_{0}^{-1}\mu(D_{l,x})$, we have
$$
(w^{1+\varepsilon})_{Q}\leq\frac{1}{\mu(Q)}
\int_{D_{l,x}}w^{1+\varepsilon}\,\mu(dx)
\leq N_{0}(w^{1+\varepsilon})_{D_{l,x}},
$$
and by \eqref{11.17.4} we have
$w_{D_{l,x}}\leq [w]_{A_{p}}N_{0}^{p-1}w_{Q}$.
This convinces us that we may concentrate on the  case
that $Q=D_{l,x}$. 

Changing scales, perhaps different
along different axes (and using the rules of changing 
the variables in the
Lebesgue integrals),
  if necessary,
convinces us that we may also assume that
$Q=D_{1,x}$. Then we observe that  $I_{Q}$ is the pointwise
limit as $\varepsilon\downarrow 0$ of the indicators
of $x+[ \varepsilon,1)^{d}$, so that we only need
to prove \eqref{11.25.4} for $Q=x+[ \varepsilon,1)^{d}$
and changing the  scales again
reduces  our task to proving \eqref{11.25.4}
only for $Q=C_{1,x}=x+[0,1)^{d}$.
In that case, as is easy to see,
 we can take $x$ as the new origin,
replace $\Omega$ with $\Omega'=[0,\infty)^{d}$, and use
the results of Section \ref{section 06.5.1.1}
in the new setting to prove \eqref{11.25.4}
for $Q=[0,1)^{d}$. 

Denote $\lambda=2N_{0}$   and for $k=0,1,...$
  define $\bar w=wI_{Q}$ and
$$
\tau_{k}(x)=\inf\{n\in\bZ:\bar w_{|n}(x)>\lambda^{k}w_{Q}
   \}.
$$
 Note that $\tau_{k+1}\geq\tau_{k}$ 
 and, obviously, $\tau_{0}=\infty$ outside $Q$.

Also observe that if $C\in \cF_{\tau_{k}}\cap \bC$,  
 then by the first inequality in
\eqref{06.5.1.3}
$$
I:=\int_{\Omega'}\bar w I_{C,\tau_{k+1}<\infty}\,\mu(dx)
\leq \int_{\Omega'}\bar w I_{C,\tau_{k }<\infty}\,\mu(dx)
\leq N_{0}\lambda^{k}w_{Q}\mu(C).
$$
On the other hand, by the second part of \eqref{06.5.1.3}
$$
I\geq  \lambda^{k+1}w_{Q}\mu(C\cap\{\tau_{k+1 }<\infty\}).
$$
Thus,
$$
\frac{\mu(C\cap\{\tau_{k+1 }<\infty\})}{\mu(C)}\leq\frac{1}{2}.
$$
By \eqref{11.24.6}
$$
w(C\cap\{\tau_{k+1 }<\infty\})\leq \beta w(C)
$$
and, since $\{\tau_{k}<\infty\}\in \cF_{\tau_{k}}$,
the set $\{\tau_{k}<\infty\}$ is represented as a
disjoint union of $C\in \cF_{\tau_{k}}\cap \bC$, so that
$$
w( \{\tau_{k+1 }<\infty\}) \leq \beta w
(\{\tau_{k }<\infty\}).
$$

It follows that  
$$
w( \{\tau_{k  }<\infty\})\leq \beta^{k}
w( \{\tau_{0}<\infty\})=\beta^{k} w(Q)
=\mu(Q)\beta^{k}w_{Q}.
$$

Now, observe that ($\mu$-a.e.) on the set
$\{\tau_{k}<\infty,\tau_{k+1}=\infty\}$
we have $w\leq \lambda^{k+1}w_{Q}$. Also
($\mu$-a.e.)
$$
\bigcup_{k=0}^{\infty}\{\tau_{k}<\infty,\tau_{k+1}=\infty\}
=\{\tau_{0}<\infty \}
$$
since $w<\infty$ ($\mu$-a.e.). Obviously, 
the terms in the above union
are disjoint. Therefore,
$$
\int_{Q}w^{1+\varepsilon}\,\mu(dx)
=\int_{Q}w^{1+\varepsilon}I_{\tau_{0}=\infty}\,\mu(dx)
+\sum_{k=0}^{\infty}\int_{Q}
w^{1+\varepsilon}I_{\tau_{k}<\infty,\tau_{k+1}=\infty}\,\mu(dx)
$$
$$
\leq \mu(Q) (w_{Q})^{1+\varepsilon}
+\sum_{k=0}^{\infty}\lambda^{(k+1) \varepsilon }
 (w_{Q})^{\varepsilon}
\int_{Q}
w I_{\tau_{k}<\infty,
\tau_{k+1}=\infty}\,\mu(dx)
$$
$$
\leq \mu(Q) (w_{Q})^{1+\varepsilon}
\Big[1+\sum_{k=0}^{\infty}\lambda^{(k+1) \varepsilon }
\beta^{k}\Big].
$$
We see that it only remains to choose $\varepsilon>0$
so small that the last series converges.
The theorem is proved.

The following somewhat sharper statement than  
\eqref{11.24.6} will be needed only
in the applications of the theory of $A_{p}$-weights
 rather than in the proof of the
Rubio de Francia
theorem.

\begin{corollary}
                                     \label{corollary 11.21.1}
If $p\in(1,\infty)$
and $w \in A_{p}$, then there exists
$\beta=\beta(p,N_{0},[w]_{p})\in (0,1)$
and $N=N(p,N_{0},[w]_{p})$ such that
for any Borel $S\subset Q\in\bQ$ we have
\begin{equation}
                                         \label{11.24.2}
\frac{w(S)}{w(Q)}\leq N\Big(\frac{\mu(S)}{\mu(Q)}\Big)^{\beta}.
\end{equation}

\end{corollary}

Indeed, by H\"older's inequality
$$
w(S)=\int_{S}w \,\mu(dx)
\leq \Big(\int_{Q}w^{1+\varepsilon}
\,\mu(dx)\Big)^{1/(1+\varepsilon)}\,
\mu^{\varepsilon/(1+\varepsilon)}(S)
$$
$$
=  \mu^{1/(1+\varepsilon)}(Q)\big[\big(w^{1+\varepsilon}\big)_{Q}\big]
^{1/(1+\varepsilon)}\mu^{\varepsilon/(1+\varepsilon)}(S)
$$
$$
\leq N\mu^{1/(1+\varepsilon)}(Q)
 w_{Q} \mu^{\varepsilon/(1+\varepsilon}(S)
=Nw(Q)\mu^{-\varepsilon/(1+\varepsilon)}(Q)
\mu^{\varepsilon/(1+\varepsilon)}(S),
$$
which is what is claimed with $\beta=\varepsilon/(1+\varepsilon)$.

For the reader's orientation we point out that
the weights satisfying \eqref{11.24.2}
are called $A_{\infty}$-weights.
 We are not going to use the remarcable
fact that $A_{\infty}=\bigcup_{p\in[1,\infty)}A_{p}$.
This fact follows from the observation that, actually,
in the proof of Theorem \ref{theorem 11.17.3}
and of \eqref{11.24.2}
only the implication
\eqref{11.24.6} was used, which turns out to be
(almost trivially) reversable.

The following result is used not only to prove
the Hardy-Littlewood weghted theorem but also
in many places in applications to partial
differential equations with almost VMO
leading coefficients.

\begin{theorem}[Self improving property]
                                 \label{theorem 11.18.3}
If $p\in(1,\infty)$
and $w \in A_{p}$, then there exists 
$$
q=q(p,N_{0},
[w]_{A_{p}})\in (1,p)
$$
 such that $w\in A_{q}$.
Furthermore, $[w]_{A_{q}}$ is estimated by a constant
$N=N(p,N_{0},
[w]_{A_{p}})$.
\end{theorem}

Proof. Note that for $p\in(1,\infty)$
the condition $w\in A_{p}$ is equivalent to
$w^{-1/(p-1)}\in A_{p/(p-1)}$ and
$$
[w^{-1/(p-1)}]_{A_{p/(p-1)}}=[w]_{A_{p}}^{1/(p-1)}.
$$
By Theorem \ref{theorem 11.17.3} 
there exist $\varepsilon>0$
and $N_{1}$ depending only on $p,N_{0},
[w]_{A_{p}}$ such that
$$
\big(w^{-(1+\varepsilon)/(p-1)}\big)_{Q}
\leq N_{1}\big(\big(w^{-1/(p-1)})_{Q}\big)^{1+\varepsilon}.
$$
Obviously, $(1+\varepsilon)/(p-1)=1/(q-1)$
for some $q\in (1,p)$ and the above inequality
means that
$$
\big(\big(w^{-1/(q-1)})_{Q}\big)^{q-1}
\leq N^{q-1}_{1}\big(\big(w^{-1/(p-1)})_{Q}\big)^{p-1}.
$$
By multiplying both parts by $w_{Q}$, we get that
$w\in A_{q}$ and $[w]_{A_{q}}\leq N^{q-1}_{1}[w]_{A_{p}}$.
The theorem is proved.

Define the Hardy-Littlewood maximal operator by
\begin{equation}
                           \label{12.28.2}
\bM f(x)=\sup_{Q\in\bQ,Q\ni x}(|f|)_{Q}=\sup_{Q\in\bQ}
I_{Q}(x)(|f|)_{Q}.
\end{equation}
Since the $Q$'s are open, $I_{Q}(x)$ are lower 
semicontinuous, so is
$\bM f$ and, hence, it is Borel measurable.
Also for an $A_{p}$-weight $w$ set
$$
\bM_{w} f(x)=\sup_{Q\in\bQ,Q\ni x}\frac{1}{w(Q)}
\int_{Q}|f|\,w(dy).
$$

\begin{lemma}
                                 \label{lemma 11.17.2}
Let $p\in[1,\infty)$ and $w \in A_{p}$. Then
for any $Q\in\bQ$ and measurable $f\geq 0$ we have
\begin{equation}
                                    \label{11.17.3}
 (f_{Q})^{p}w_{Q}\leq [w]_{A_{p}}(f^{p}w)_{Q}.
\end{equation}
\end{lemma}

Proof. If  $p=1$, \eqref{11.17.3} follows from 
\eqref{11.17.4} and the linearity of the integral.
If $p>1$, by H\"older's inequality
$$
(f_{Q})^{p}w_{Q}=\big((fw^{1/p})w^{-1/p}\big)_{Q}^{p}w_{Q}
\leq (f^{p}w)_{Q}\Big(\big(w^{-1/(p-1)}\big)_{Q}\Big)^{p-1}
w_{Q}
$$
and \eqref{11.17.3} follows in light of our definitions.
The lemma is proved.

Observe that \eqref{11.17.3} implies the following.
\begin{corollary}
                                     \label{corollary 11.21.4}
For $p\in[1,\infty)$, $w\in A_{p}$, and Borel $f\geq0$
we have
$$
(\bM f)^{p}\leq[w]_{A_{p}}\bM_{w}( f ^{p}).
$$  
\end{corollary}

\mysection{The case of $\mu$
with a ``real'' doubling property on $\bQ$}
                                 \label{section 11.21.1}
In this section we assume that $\mu$
satisfies the conditions
\eqref{11.28.1} and \eqref{11.24.4} from
 Sections \ref{section 06.5.1.1} 
and \ref{section 11.20.1} 
and 
satisfies an even stronger doubling condition
than \eqref{11.24.1} 
because we are going to use
a covering lemma similar to the Vitali one.
Let $\Omega^{o}$ be the interior of $\Omega$.

Introduce $\bQ$
as the collection of
$$
C_{l}(x)= x+(-l^{k_{1}}/2,l^{k_{1}}/2)\times...\times
(-l^{k_{d}}/2,l^{k_{d}}/2)  ,\quad x\in\Omega,l>0,
$$
that lie entirely in $\Omega$. Actually, this is the same
collection as the collection of $D_{l,x},x\in\Omega,l>0$,
from Section \ref{section 11.20.1}. We use a different 
notation which allows us to express the doubling
condition more easily.
Suppose that for any $C_{l}(x)\in \bQ$ 
we have 
\begin{equation}
                                    \label{11.21.3}
\mu( C_{4l}(x)\cap\Omega  )\leq N_{0}\mu( C_{l}(x) ).
\end{equation}
All $A_{p}$-weights in this section are
$A_{p}$-weights relative to these $\bQ$, $\mu$.

Observe that, if $C_{l}(x)\in \bQ$, then $C_{l}(x)=D_{l,y}$
for a $y\in\Omega$ and $D_{2l,y}\subset C_{4l}(x)\cap\Omega^{o}$.
Therefore, condition \eqref{11.21.3} implies that 
 condition \eqref{11.24.1} is satisfied
by the above $\bQ$ and we can use the results of Section
\ref{section 11.20.1}. In particular,
\eqref{11.17.4} implies the following.

\begin{corollary}
                                     \label{corollary 11.17.1} 
For any $p\in[1,\infty)$ and $w\in A_{p}$,
the measure $w(A) $ satisfies
the doubling condition  \eqref{11.21.3}
 with constant 
$N_{0}^{p}[w]_{A_{p}}$.  
\end{corollary}

\begin{lemma}[Maximal inequality]
                                 \label{lemma 11.21.3}
If $p\in[1,\infty)$ and $w \in A_{p}$, then
there exists a constant $N$, depending only
on $p,N_{0}$, and $[w]_{A_{p}}$,
such that for any $\lambda>0$ and 
Borel $f$
\begin{equation}
                                    \label{11.17.5}
 w(\{x:\bM_{w} f(x)> \lambda\})\leq N\lambda^{-1}
\int_{\Omega}|f|\,w(dx).
\end{equation}
\end{lemma}

Proof.  
We may assume that $f\geq0$
and $f\in L_{1}(w)$. Define $A(\lambda)
=\{x\in\Omega:\bM_{w} f(x)>\lambda\}$. This set is open
since $\bM_{w} f$ is lower semicontunuous.
  Take a compact set $K\subset
A(\lambda)$. Then by the definition of $A(\lambda)$ for
any $x\in K$ there exists a $Q\in\bQ$ such that
$x\in Q$ and $\int_{Q}f\,w(dx)>\lambda
w(Q)$. Then, of course,
$Q\subset A$. 

By the compactness of $K$, there is a finite collection
$Q_{1},...,Q_{n}\in\bQ$ covering $K$
and such that for each $Q=Q_{i}$
we have $\int_{Q}f\,w(dx)>\lambda
w(Q)$.

Now we use a Vitali's  covering argument. If $Q
=C_{l}(x)\in\bQ$, then
define its size to be $l$ and set
$ Q'=C_{4l}(x)\cap\Omega^{o}$.   Then denote by $\hat Q_{1}$
any of $Q_{i}$ which has the largest size
 and set it aside.
  Next,
introduce $\hat Q_{2}$ as one of the remaining $Q_{i}$
which has the  largest size
 between those $Q_{i}$ that have
{\em no\/} intersection with $\hat Q_{1}$. It may happen that
no such $Q_{i}$ exist. Then it is almost obvious that
$Q_{i}\subset\hat Q'_{1}$ for any $i$.
If $\hat Q_{2}$ exists, we proceed further.

If we have already defined $\hat Q_{1},...,\hat Q_{k}$,
then we define $\hat Q_{k+1}$ as one of the ``cubes'' in
the family  of
 \begin{equation}
                                            \label{06.10.18.3}
\{Q_{1},...,Q_{n}\}
\setminus\{\hat Q_{1},...,\hat Q_{k}\},
\end{equation} 
 which are disjoint from $\hat Q_{1},...,\hat Q_{k}$, which has the
largest size (between those in \eqref{06.10.18.3}, that 
are  disjoint from $\hat Q_{1},...,\hat Q_{k}$).
In finitely many steps we will come to a $k$
 for which
any  ``cube'' in family \eqref{06.10.18.3} intersects
one of $\hat Q_{1},...,\hat Q_{k}$ or else
the family is empty.
In the second case, obviously, for any $i$
  \begin{equation}
                                            \label{06.10.18.4}
Q_{i}\subset\bigcup_{j=1}^{k}\hat Q_{j}'.
\end{equation} 

It turns out that \eqref{06.10.18.4} also holds for any $i$
in the first case. Indeed, if, for a fixed $i$,
 $Q_{i}$ has a nonempty intersection
with a $\hat Q_{j}$, then define $r=r(i)$ as the smallest such $j$
and observe that, if $r=1$, then as has been pointed out above,
$Q_{i}\subset \hat Q_{1}'$ and \eqref{06.10.18.4} holds.
If $r>1$, then the size of
$ Q_{i}$ is no greater than that of $\hat Q_{r}$
by the choice of $\hat Q_{r}$ and because $Q_{i}$
has no intersection with $\hat Q_{1},...,\hat Q_{r-1}$
by the definition of $r$. Now as above,
this combined with $\hat Q_{r}\cap Q_{i}\ne\emptyset$
implies that $Q_{i}\subset \hat Q_{r}'$. This proves
\eqref{06.10.18.4}.

  It follows that
 $$
K\subset\bigcup_{j=1}^{k}\hat Q_{j}'.
 $$
Finally, use the fact that, owing to the doubling property
in Corollary \ref{corollary 11.17.1}, we have
 $
w(Q')\leq Nw(Q),
$
where the constant $N$ depends 
only on $p,N_{0}$, and $[w]_{A_{p}}$,
so that
$$
w(K)\leq\sum_{j=1}^{k}w(\hat Q_{j}')
\leq N\sum_{j=1}^{k}w(\hat Q_{j} )
$$
$$
\leq N \lambda^{-1}\sum_{j=1}^{k}
\int_{\hat Q_{j}}f\,w(dx)
\leq N \lambda^{-1}\int_{\Omega}f\,w(dx).
$$

By taking a sequence of compact sets $K_{m}\uparrow
A(\lambda)$ and passing to the limit, we get 
\eqref{11.17.5} in our particular case.
 The lemma is proved.

We can now prove one of the fundamental results of the theory.

\begin{theorem}
                                \label{theorem 11.21.3}
(i) If $p\in[1,\infty)$ and $w\in A_{p}$, then the operator
$\bM f$ is of weak $(p,p)$-type with respect
to measure $w$. More precisely,  there exists
a constant $N$ depending only on $p,N_{0}$, and 
$[w]_{A_{p}}$ such that for any $\lambda>0$ and Borel 
$f\geq0$ we have
$$
w(\{x:\bM f>\lambda\}\leq N\lambda^{-p}
\int_{\Omega}f^{p}\,w(dx).
$$

(ii) If $p\in(1,\infty)$ and $w\in A_{p}$, then  
$\bM f$ is  a bounded operator in $L_{p}(w)$,
 that is there exists
a constant $N$ depending only on $p,N_{0}$, and 
$[w]_{A_{p}}$ such that for any   
$f\in L_{p}(w)$ we have
$$
\|\bM f\|_{L_{p}(w)}\leq N\|  f\|_{L_{p}(w)}.
$$

\end{theorem}

Proof. (i) By Lemma \ref{lemma 11.21.3} and
Corollary \ref{corollary 11.21.4}
$$
w(\{x:\bM f>\lambda^{1/p}\}
\leq w(\{x:\bM_{w}( |f|^{p})>[w]_{A_{p}}^{-1}\lambda  \}
\leq N[w]_{A_{p}}\lambda^{-1}
\int_{\Omega}f^{p}\,w(dx)
$$
and it only remains to make an obvious substitution.

(ii) By Theorem \ref{theorem 11.18.3}
 there exists $q=q(p,N_{0},
[w]_{A_{p}})\in (1,p)$ such that $w\in A_{q}$.
Then assertion (i) implies that 
$\bM f$ is a weak $(q,q)$-operator 
with respect to   measure $w$. The fact that it also maps
bounded functions into bounded ones with norm
one is almost obvious. Now our assertion follows 
from the Marcinkiewicz
interpolation theorem. The theorem is proved.

It is worth noting that in \cite{Le_08}
one can find an elegant and short proof of Theorem
\ref{theorem 11.21.3} (ii)
not using Theorem \ref{theorem 11.18.3}. 

\mysection{Rubio de Francia extrapolation theorem}
                                       \label{section 11.21.10}
We work in the framework of Section 
\ref{section 11.21.1} and suppose that $\mu$
satisfies \eqref{11.28.1}, \eqref{11.24.4}, 
 and \eqref{11.21.3}.

\begin{theorem}
                                 \label{theorem 11.21.6}
Let $p,q\in(1,\infty)$, $ K_{q}$ be
a constant  from $(1,\infty)$ and $f,g$ be Borel
nonnegative functions of $\Omega$.
Then there exists $K_{p}=K (p,q,K_{q}, N_{0})  
\geq 1$ such that if
\begin{equation}
                                     \label{11.21.6}
\int_{\Omega}f^{p} w_{p}\,\mu(dx)
\leq \int_{\Omega}g^{p} w_{p}\,\mu(dx)
\end{equation}
for any $w_{p}\in A_{p}$ with $[w]_{A_{p}}\leq K_{p}$,
then
\begin{equation}
                                     \label{11.21.7}
\int_{\Omega}f^{q} w_{q}\,\mu(dx)
\leq 4^{q}\int_{\Omega}g^{q} w_{q}\,\mu(dx)
\end{equation}
for any $w_{q}\in A_{q}$ with $[w]_{A_{q}}\leq K_{q}$.
\end{theorem}

Proof.
 We follow the proof of Theorem 2.5
of \cite{DK_18} or the proof of Theorem 1.4
of \cite{C-UMP_11}, which are streamlined versions
of the proof of Theorem IV.5.19 of \cite{GR_85}.
 The proof is  striking.  
A slightly different proof one can find in 
\cite{Du_11}. See also \cite{Sa_15} for
generalizations in the case of abelian groups. 

Denote by $N_{1}$ the constant $N$ which suits
Theorem \ref{theorem 11.21.3} (ii)
with $q$ in place of $p$ for any $w\in A_{q}$
with $[w]_{A_{q}}\leq K_{q}$. Also set
$q'=q/(q-1)$ and denote by $N_{2}$
the constant $N$ which suits
Theorem \ref{theorem 11.21.3} (ii)
with $q'$ in place of $p$ for any $w\in A_{q'}$
with $[w]_{A_{q'}}\leq K^{1/(q-1)}_{q}$.
Set
$$
K_{p}=2^{p}N_{1}^{p-1}N_{2},
$$
and assuming that
\eqref{11.21.6} holds for all $w\in A_{p}$
with $[w]_{A_{p}}\leq K_{p}$
 we prove \eqref{11.21.7} for any $w_{q}\in A_{q}$
such that $[w_{q}]_{A_{q}}\leq K_{q}$.

We fix such a $w_{q}$ and for 
nonnegative
$h\in L_{q}(w_{q})$  
  define
$$
\bR h(x)=\sum_{k=0}^{\infty}\frac{\bM^{k}h(x)}
{2^{k}N_{1}^{k}},
$$
where $\bM^{k}$ is the $k$-th iteration of $\bM$:
$\bM^{0} h=h$, $\bM^{k+1}=\bM\bM^{k}$.

Observe that, obviously, $h\leq \bR h$,
$\|\bR h\|_{L_{q}(w_{q})}\leq 2\|  h\|_{L_{q}(w_{q})}$,
$\bM\bR h\leq 2N_{1} \bR h$, so that,
  if   $\|h\|_{L_{q}(w_{q})}>0$, then $\bR h\in A_{1}$
and
$$
[\bR h]_{A_{1}}\leq 2N_{1}.
$$

Next,  we know that
$w_{q}^{1-q'}\in A_{q'}$  and
 $\bM$ is bounded on $L_{q'}(w_{q}^{1-q'})$. 
Since
$[w_{q}^{1-q'}]_{A_{q'}}=([w_{q}]_{A_{q}})^{1/(q-1)}\leq 
K^{1/(q-1)}_{q}$, we have that 
$$
\|\bM h\|_{L_{q'}(w_{q}^{1-q'})}\leq
N_{2}\| h\|_{L_{q'}(w^{1-q'})}.
$$
This means that the operator $\bM'h=w_{q}^{-1}\bM(hw_{q})$
is bounded in $L_{q'}(w_{q} )$ with norm less than $N_{2}$. Indeed,
 for $h\in L_{q'}(w_{q})$ we have
$$
\|\bM'h\|_{L_{q'}(w_{q})}=\|\bM (hw_{q})\|_{L_{q'}(w_{q}^{1-q'})}
\leq N_{2}\|hw_{q}\|_{L_{q'}(w_{q}^{1-q'})}
=N_{2}\| h\|_{L_{q'}(w_{q})}.
$$

This allows us to introduce the operator
$$
\bR' h(x)=\sum_{k=0}^{\infty}\frac{(\bM')^{k}h(x)}
{2^{k}N_{2}^{k}},\quad h\in L_{q'}(w_{q}), h\geq0,
$$
and claim that $h\leq \bR' h$,
$\|\bR' h\|_{L_{q'}(w_{q})}\leq 2\|\bR' h\|_{L_{q'}(w_{q})}$,
$\bM'\bR' h\leq 2N_{2} \bR' h$, so that,
if $\|h\|_{L_{q'}(w_{q})}>0$, then $w_{q}\bR' h\in A_{1}$
and
$$
[w_{q}\bR' h]_{A_{1}}\leq 2N_{2}.
$$

Then observe that wile proving \eqref{11.21.7}
we
  may certainly assume that $0<\|g\|_{L_{q}(w_{q})}<\infty$.
Take a nonnegative test function $h\in L_{q'}(w_{q})$
such that $\|h\|_{L_{q'}(w_{q})}>0$.
From the properties of $\bR$ and $\bR'$ 
and Lemma \ref{lemma 11.19.1} we get that
the function
$$
 w_{p}:=(\bR g)^{1-p}w_{q}\bR'h
$$
is an $A_{p}$-weight and
$$
[ w_{p}]_{A_{p}}\leq [\bR g]_{A_{1}}^{p-1}[w_{q}\bR'h]_{A_{1}}
\leq 2^{p}N_{1}^{p-1}N_{2}=K_{p}.
$$
 
Furthermore, since $\|\bR g\|_{L_{q}(w_{q})}\leq 2
\|  g\|_{L_{q}(w_{q})}$ and $\|\bR' h\|_{L_{q'}(w_{q})}\leq 2
\|  h\|_{L_{q'}(w_{q})}$, by H\"older's inequality
\begin{equation}
                                     \label{11.21.8}
\int_{\Omega}(\bR g)w_{q}\bR'h \,\mu(dx)
\leq 4\|  g\|_{L_{q}(w_{q})}\|  h\|_{L_{q'}(w_{q})}.
\end{equation}

Now, by H\"older's inequality and by assumption
$$
\int_{\Omega}fhw_{q}\,\mu(dx)
=\int_{\Omega}f(hw_{q}/ w_{p})\, w_{p}(dx)
\leq\|f\|_{L_{p}( w_{p})}\|hw_{q}/ w_{p}
\|_{L_{p/(p-1)}( w_{p})}
$$
$$
\leq\|g\|_{L_{p}( w_{p})}\|hw_{q}/ w_{p}
\|_{L_{p/(p-1)}( w_{p})},
$$
where, in light of $g\leq \bR g$ and \eqref{11.21.8},
$$
\|g\|^{p}_{L_{p}( w_{p})}
=\int_{\Omega}g^{p}(\bR g)^{1-p} \bR'h\,w_{q}(dx)
$$
$$
\leq \int_{\Omega} (\bR g)  \bR'h\,w_{q}(dx)
\leq 4\|  g\|_{L_{q}(w_{q})}\|  h\|_{L_{q'}(w_{q})},
$$
and, in light of $h\leq \bR' h$ and \eqref{11.21.8},
$$
\|hw_{q}/ w_{p}\|^{p/(p-1)}_{L_{p/(p-1)}( w_{p})}=
\int_{\Omega}h^{p/(p-1)} (\bR g)^{p }(\bR' h)^{p/(1-p)}
(\bR g)^{1-p} \bR'h\,w_{q}(dx)
$$
$$
\leq \int_{\Omega} (\bR g)  \bR'h\,w_{q}(dx)
\leq 4\|  g\|_{L_{q}(w_{q})}\|  h\|_{L_{q'}(w_{q})}.
$$

Hence,
$$
\int_{\Omega}fh \,w_{q}(dx)\leq
4\|  g\|_{L_{q}(w_{q})}\|  h\|_{L_{q'}(w_{q})}
$$
and the arbitrariness of $h\in L_{q'}(w)$
proves \eqref{11.21.7}. The theorem
is proved.

To extract the most important for us
consequences of this theorem we split
the coordinates of $x$ into 
several groups.

We take integers $m\geq2, d_{1},...,d_{m}\geq1$ such that
$d_{1}+...+d_{m}=d$, define $l_{0}=0$, $l_{i+1}=l_{i}+d_{i+1}$,
$i=0,...,m-1$,  and 
express the points in $\bR^{d}$ as
$$
x=(x_{1},\ldots,x_{d})=(\check x_{1},\ldots,\check x_{m}),
$$
where $\check x_{i}=(x_{l_{i-1}+1},\ldots, x_{l_{i}})$.
Then $\bR^{d}=\bR^{d_{1}}\times...\times \bR^{d_{m}} $
and, accordingly,
$\Omega=\Omega_{1}\times...\times \Omega_{m} $,
where $\Omega_{i}$ is the projection of $\Omega$
on $\bR^{d_{i}}$. We use the notation 
$\check x_{i}=(x_{l_{i-1}+1},\ldots, x_{l_{i}})$
for generic points in $\bR^{d_{i}}$
and set
$$
\hat \Omega^{i}=
\Omega_{i}\times\cdots\times \Omega_{m}
=\{\hat x^{i}=(\check x_{i},...,\check x_{m})\}.
$$

Let $\bQ_{i}$ be the family of
 projections on $\bR^{d_{i}}$
of elements of $\bQ$ from Section \ref{section 11.21.1},
so that $\bQ=\bQ_{1}\times...\times \bQ_{m}$ and if
$C_{l}(x)\in\bQ$, then $C_{l}(x)=C_{l}(\check x_{1})
\times...\times C_{l}(\check x_{m})$, where
$C_{l}(\check x_{i})\in\bQ_{i}$. We set
$$
\bQ^{i}=\bQ_{i}\times...\times\bQ_{m}.
$$
  We assume that, for each $i$, we are given
a measure $\mu_{i}$, which satisfies conditions
\eqref{11.28.1} and  \eqref{11.24.4}  relative to
  $\Omega_{i}$ and satisfy the following
doubling condition: if
$C_{l}(\check x_{i})\in\bQ_{i}$, then
\begin{equation}
                                            \label{11.27.1}
\mu_{i}(C_{4l}(\check x_{i})\cap\Omega_{i})
\leq N_{0}^{1/m}\mu_{i}(C_{ l}(\check x_{i})).
\end{equation}
Observe that $N_{0}$ in \eqref{11.24.1}
or $N_{0}^{1/m}$ in \eqref{11.27.1}
are not necessarily the best constants
for which these conditions are valid and we use
$N_{0}^{1/m}$ in \eqref{11.27.1} just in order to make
our last assumption that
 $\mu=\mu_{1}\times...\times \mu_{m}$
somewhat consistent. We set
$$
\mu^{i}=\mu_{i}\times...\times\mu_{m}.
$$

 Next, we introduce $A_{p}(\bQ_{i},\mu_{i})$
as $A_{p}$-weights relative 
to $\bQ_{i},\mu_{i}$. One easily checks that if $w_{i}$
are $A_{p}$-weights relative to $\bQ_{i},\mu_{i}$, then
$w_{i}w_{j}$ is an  $A_{p}$-weights relative to 
$\bQ_{i}\times\bQ_{j},\mu_{i}\times\mu_{j}$ and
$$
[w_{i}w_{j}]_{A_{p}}\leq[w_{i} ]_{A_{p}}[ w_{j}]_{A_{p}}.
$$
 
Finally, for $i\in\{1,...,m\}$,
 $p_{i}\in(1,\infty)$, and
weights $w_{i}$ given on $\Omega_{i}$ introduce
$$
\| f
\|_{L_{p_{1},\ldots,p_{m}}(w_{1},\ldots,w_{m})}^{p_{m}}
$$
$$
:=\int_{\Omega_{m}}\Big(\cdots
\Big(\int_{\Omega_{2}}
\Big(\int_{\Omega_{1} } |f|^{p_{1}}w_{1}\,\mu_{1}(d\check x_{1})
\Big)^{p_{2} /p_1 }w_{2}\,\mu_{2}(d\check x_{2})\Big)^{p_{3}/p_{2}}
$$
\begin{equation}
                                       \label{4.9.60}
\cdots\Big)^{p_{m} /p_{m-1} }w_{m}\,\mu_{m}(d\check x_{m}).
\end{equation}

 \begin{theorem}
                           \label{theorem 12.2.1}
Let $K^{*},p_{k}\in(1,\infty)$, $w_k\in
A_{p_{k}}( \Omega_k,\mu_{k})$,
 $[w_k]_{p_{k}}\le K^{*}$, $k=1,\ldots,m$, and
let
$u,g$ be measurable functions on $\Omega$. Then there exists
a constant $K_0=K_0(d,m,p_1,\ldots,p_{m},
N_{0},K^{*})\ge
1$ such that if
$$
\|u\|_{L_{p_{1}}(w )}
\le  \|g\|_{L_{p_{1}}(w )}
$$  for every $w\in
A_{p_{1}}(\bQ,\mu)$ with
 $[w]_{p_{1}}\le K_0$, then we have
$$
\|u\|_{L_{p_{1},\ldots,p_{m}}(w_1,\ldots,w_{m})}\le
4^{m}\|g\|_{L_{p_{1},\ldots,p_{m}}(w_1,\ldots,w_m)}.
$$  
\end{theorem}

Proof. We follow the proof of Corollary 2.7
in \cite{DK_18} or Theorem 8.1 in \cite{DKr_18}.  
Take the function $K (p,q,K_{q}, N_{0})$
from Theorem \ref{theorem 11.21.6}
and restate it in the following way.

If for fixed $0\leq j \leq m-1$, $K_{j+1}\geq1$, 
 and two nonnegative functions
$U_{j}$ and $G_{j}$ on $\hat\Omega^{j+1}$
 it holds that (a) 

\begin{equation}
                                         \label{4.22.3}
\int_{\hat\Omega^{j+1} }
U_{j}^{p_{j}}w(\hat x^{j+1} )\,
 \mu^{j+1}(d\hat x^{j+1})   \leq 4^{(j-1)p_{j}} \int_{\hat\Omega^{j+1} }
G_{j}^{p_{j}}w(\hat x^{j+1} )\,
\mu^{j+1}(d\hat x^{j+1})
\end{equation}
   for every $w\in
A_{p_{j}}(\bQ^{j+1},\mu^{j+1} )$ as long as
\begin{equation}
                                         \label{11.27.3}
[w]_{p_{j}}\le K_{j}=K(p_{j},p_{j+1},K^{*}K_{j+1},N_{0}),
\end{equation}
then we have (b) 
\begin{equation}
                                         \label{4.22.4}
\int_{\hat\Omega^{j+1} }
U_{j}^{p_{j+1}}w(\hat x^{j+1} )\,
\mu^{j+1}(d\hat x^{j+1})  \leq  4 ^{jp_{j+1}}\int_{\hat\Omega^{j+1} }
G_{j}^{p_{j+1}}w(\hat x^{j+1} )\,
\mu^{j+1}(d\hat x^{j+1})
\end{equation}
 for every $w\in
A_{p_{j+1}}(\bQ^{j+1},\mu^{j+1} )$ with
 $[w]_{p_{j+1}}\le K^{*}K_{j+1}$.

We define  $K_{m }=1$  and define
  $K_{j}$ for $j=0,1,...,m-1$ recurrently by
the equation in \eqref{11.27.3}.

Also set  $U_{0}(x)=u(x)$,
$$
U_{j}(\hat x^{j+1})=
\Big(\int_{ \Omega_{j}}U^{p_{j}}_{j-1}(\hat x^{j})
\, w_{j}(\check x_{j})\,\mu_{j}(d\check x_{j})\Big)^{1/p_{j}},\quad  1\leq j\leq m-1,
$$
and similarly we introduce
$G_{j}$'s
by taking $g$ in place of $u$. To prove the theorem,
it suffices to prove that (b) holds for
$j=m-1$ because $w_{m}\in A_{p_{m}}(  \bQ_{m},\mu_{m})$
and $[w_{ m}]_{A_{p_{m}}}\leq K^{*}=K^{*}K_{m}$.
 We are going to use the induction on $j=0,1,\ldots,m-1$.

Observe that
 (b) holds for $j=0$   by assumption.
Suppose that it holds for a $j\in \{0,1,\ldots,m-2\}$.
Then \eqref{4.22.4} also holds for
$$
w(\hat x^{j+1} ):=w_{j+1}(\check x_{j+1})w(\hat x^{j+2} )
$$
if $w(\hat x^{j+2} ) \in A_{p_{j+1}}(\bQ^{j+2},\mu^{j+2})$ and $[w(\hat x^{j+2} 
)]_{A_{p_{j+1}}}\le K_{j+1} $,
because then
$
[w(\hat x^{j+1} )]_{A_{p_{j+1}}}\le K^{*}K_{j+1}
 $. Remarkably, this implies that
(a) holds with $j+1$ in place of $j$.
Then (b) also holds with $j+1$ in place of $j$.
This justifies the induction and proves the theorem.

\begin{remark}
                                 \label{remark 12.21.1}
By relabeling the coordinates one sees
that Theorem \ref{theorem 12.2.1} holds true
if the repeated integrals in \eqref{4.9.60}
are taken in any other order.
\end{remark}

\mysection{Generalized Fefferman-Stein theorem}
                                       \label{section 11.29.1}

We keep working in the setting of Section \ref{section 11.21.1}  
with $\mu$ satisfying  \eqref{11.28.1},
  \eqref{11.24.4}, and \eqref{11.21.3}. All
$A_{p}$-weights are also taken from there.
Since, obviously, $\cM f\leq\bM f$
(cf.~\eqref{12.28.1} and \eqref{12.28.2})
the following is a 
corollary of Theorem \ref{theorem 11.21.3} (ii).

\begin{theorem}[Hardy-Littlewood]
                                           \label{theorem 12.21.2}
Let $p\in(1,\infty)$ and $w\in A_{p}$. Then      
for any $f\in \cL_{p}(w)$ we have  
\begin{equation}
                                      \label{12.19.4}
\|\cM f\|_{\cL_{p}(w)}
\leq N\|f \|_{\cL_{p}(w)},
\end{equation}
where $N$ depends only on $p$, $N_{0}$, and $[w]_{A_{p}}$. 
\end{theorem}

Next result is   Lemma 2.8 of \cite{DKr_18}.
Its proof given for the sake of
completeness is also taken from there.  
Recall that $\bC$ is introduced  
in Section~\ref{section 06.5.1.1}.

\begin{lemma}
                                         \label{lemma 12.16.1}
Let $g$ be a measurable function on
$\Omega$, $\gamma\in(0,1]$, $vI_{C}\in L_{1}(\mu)$
for any   $C\in\bC $,
and let $v_{|n}\to0$ as $n\to-\infty$ on $\Omega$. Assume that
  $|u|\leq v$
and
for any   $C\in\bC $ there exists a
measurable function  $u^{C}$ given on $C$ such that
$|u|\leq u^{C}\leq v$ on $C$ and, for any $x\in C$
\begin{equation}
                                               \label{6.29.3}
\Big( \dashint_{C}
\dashint_{C}\big|u^{C}(z)-u^{C}(y)\big|^{\gamma}
\,\mu(dz)\mu(dy)\Big)^{1/\gamma} \leq
 g (x).
\end{equation}
Let $p\in[1,\infty)$ and $w$ be an $A_{p}$-weight. 
Then for any $\lambda>0$ we have
\begin{equation}
                                              \label{6.29.1}
 w( \{x:|u(x)|\geq\lambda\}) \leq N
\nu^{-\beta}\lambda^{-\gamma\beta}
\int_{\Omega}g^{\gamma\beta}(x)I_{\cM v(x)>\alpha \lambda}\,w(dx),
\end{equation}
where   $\alpha=(2N_{0})^{-1}$ and  $\nu= 1-2^{-\gamma}$
and the constants $\beta\in(0,1)$ and  $N$ depend
 only on $p, N_{0}$, and
$[w]_{A_{p}}$. 
\end{lemma}

Proof. Obviously we may assume that   $u\geq0$.
Fix a $\lambda>0$ and define  \vspace{5pt}
$$
\tau(x)=\inf\big\{n\in\bZ:v_{|n}(x)>\alpha\lambda\big\}.
\vspace{5pt}$$
We know that $\tau$ is a stopping time and if $\tau(x)<\infty$,
then
$$
  v_{|n}(x)\leq \lambda/2,\quad\forall n\leq\tau(x).
$$
We also know that $v_{|n}\to v\geq u$ (a.e.) as $n\to\infty$. 
It follows  that (a.e.)
\begin{align*} 
&\big\{x:u(x)\geq\lambda\big\}=\big\{x:u(x)\geq\lambda,\tau(x)<\infty\big\}
\vspace{5pt}\\
&=\big\{x:u(x)\geq\lambda,  v_{|\tau}(x)\leq \lambda/2\big\}
=\bigcup_{n\in\bZ}\bigcup_{C\in
\cF^{\tau}_{n}}A_{n}(C),
\end{align*}
where
$$
A_{n}(C):=\big\{x\in C:u(x)\geq\lambda, v_{|n}(x)\leq \lambda/2\big\},
\vspace{5pt}$$
and $\cF^{\tau}_{n}$ is the family of disjoint elements
of $\bC_{n}$ such that
$$
\big\{x:\tau(x)=n\big\}=\bigcup_{C\in \cF^{\tau}_{n}}C.
$$

Next, for each $n\in\bZ$ and $C\in\bC_{n}$ on the
set $A_{n}(C)$,
if it is not empty,
we have $v_{|n}=v_{C}$ and on $A_{n}(C)$
$$
u^{\gamma}-(v_{C})^{\gamma}\geq\lambda^{\gamma}(1-2^{-\gamma})
=\nu\lambda^{\gamma} .
 $$
  We use this
and the inequality $|a-b|^{\gamma}\geq|a|^{\gamma}-|b|^{\gamma}$
and conclude that for $x\in A_{n}(C)$,
$$
\dashint_{C}\big|u^{C}(x)-u^{C}(y)\big|^{\gamma}\,\mu(dy)
\geq\big(u^{C}(x)\big)^{\gamma}-\dashint_{C}\big(u^{C}(y)\big)^{\gamma}\,\mu(dy)
 $$
$$
\geq u^{\gamma}(x)-\dashint_{C}v^{\gamma}(y)\,\mu(dy)
\geq u^{\gamma}(x)-\big(v_{C}(x)\big)^{\gamma}
\geq\nu\lambda^{\gamma},
\vspace{5pt}$$
 so that
by Chebyshev's inequality  \vspace{5pt}
$$
\mu(A_{n}(C))\leq
\nu^{-1}\lambda^{-\gamma}\int_{ C }
\dashint_{C}\big|u^{C}(z)-u^{C}(y)\big|^{\gamma}\,\mu(dz)\mu( dy ).
\vspace{5pt}$$

It follows by assumption \eqref{6.29.3} that \vspace{5pt}
$$
\frac{\mu(A_{n}(C))}{\mu( C) }\leq
\nu^{-1}\lambda^{-\gamma}g^{\gamma}(x)
\vspace{5pt}$$
for any $x\in\Omega$.     Corollary 
\ref{corollary 11.21.1} implies that
$$
w(A_{n}(C))\leq N_{w,\beta}\nu^{-\beta}\lambda^{-\gamma\beta}
g^{\gamma\beta}(x)w(C).
$$
Since this holds for any $x\in C$,
$$
w(A_{n}(C))\leq N_{w,\beta}\nu^{-\beta}\lambda^{-\gamma\beta}\int_{C}
g^{\gamma\beta}(x)\,w(dx).
$$
Hence,
\begin{align*}
&w\big\{x:u(x)\geq\lambda\big\}\leq N_{w,\beta}
\nu^{-\beta}\lambda^{-\gamma\beta}
\sum_{n\in\bZ}\sum_{C\in \cF^{\tau}_{n}}\int_{C}
g ^{\gamma\beta}\,w(dx)
\vspace{5pt}\\
&=N_{w,\beta}\nu^{-\beta}\lambda^{-\gamma\beta}
\int_{\Omega}g^{\gamma\beta}I_{\tau<\infty}\,
w(dx).
\vspace{5pt}
\end{align*}
It only remains to observe that $\{\tau<\infty\}=\{\cM
v>\alpha\lambda\}$. The lemma is proved.

\begin{remark}
                                        \label{remark 11.29.1}
It is worth saying a few words about
the history of Lemma \ref{lemma 12.16.1},
which is the 
core of the approach based on the Fefferman-Stein theorem.
In case $\gamma=1$, $u=u^{C}=v$, $w(x)\equiv1$,
and slightly more general right-hand side
of \eqref{6.29.3}, estimate \eqref{6.29.1} with $\beta=1$
becomes Lemma 3.2.9 of \cite{Kr_08} and serves there
as one of the mains tools of obtaining  a priori estimates
for Sobolev solutions of {\em linear\/} elliptic and parabolic
equations with continuous or VMO or else
 almost VMO coefficients.
The point is that \eqref{6.29.3} for $D^{2}u$ in place of 
$u^{C}$ amounts to pointwise estimating the 
so-called sharp function of
$D^{2}u$ (cf.~Theorems \ref{theorem 11.30.1}, \ref{theorem 12.1.1})
and \eqref{6.29.1}    easily leads to estimates
of the $L_{p}$-norms of $D^{2}u$ through 
the norm of its sharp function (cf.~the proof
of Theorem
\ref{theorem 12.19.1}).
In this form  Lemma \ref{lemma 12.16.1} was designed
to use the Fefferman-Stein theorem
instead of explicit singular integral
representation of the derivatives of solutions
 and subsequent application of
the Coifman-Rochberg-Weiss commutator theorem
 for singular integrals  in order 
 to treat equations with main VMO
coefficients, which was first done   
in \cite{CFL_91}   and
\cite{CFL_93}  and later continued in
an avalanche of papers. Using the Fefferman-Stein theorem
does not require any explicit representation
formulas and turned out to be applicable to
many linear and {\em fully nonlinear\/} equations.

In its initial form Lemma \ref{lemma 12.16.1}
turned out to provide a crucial information
even in the case when the coefficients
are only measurable with respect to
one variable and almost VMO with respect to
the others. This line of research started
in \cite{KK_07_1} and \cite{KK_07_2}
and continued in many papers, see, for
instance, \cite{DK_11_1,DK_11_2,DK_11_3,DK_15},
and the references therein. In particular, in
\cite{Ki_10} it is allowed for all main coefficients
to depend in a measurable  way  on time and one space
coordinate apart from one
which is supposed to be measurable in time
and VMO in spatial variables. In the elliptic case in
\cite{DKr_10}   the coefficients are measurable with respect
to two spatial variables.

The need to generalize Lemma \ref{lemma 12.16.1}
and add $u^{C}$ still with $\gamma=\beta=1$, $w\equiv1$
comes when one wants to allow the direction
in which the coefficients are only measurable to 
depend on $x$. This was first noted in
\cite{Kr_09} for nondivergence type equations
and first used in \cite{Do_10} for divergence equations.

In a remarkable paper \cite{DK_18} the authors
proved the version of Lemma \ref{lemma 12.16.1}
with $\gamma=1$ but with $A_{p}$-weights.
This allowed the authors to get mixed norms estimates
just by referring to Theorem \ref{theorem 12.2.1}.

The necessity to have $\gamma\in(0,1]$
(actually, very small one) came when we
started applying the same methodology 
to fully nonlinear equations. Then 
Lemma \ref{lemma 12.16.1} with 
$\beta=1$, $w\equiv1$, appeared in
\cite{Kr_10}.   

In its final form Lemma \ref{lemma 12.16.1}
was needed in \cite{DKr_18} when
we investigated fully nonlinear equations
in the spaces with mixed norms.

\end{remark}

Whenever it makes sense define {\em the sharp ``dyadic'' 
function of\/} $f$
 by
$$
f^{\#}_{\gamma }(x)= 
\sup_{\substack{C\in \bC ,\\
C\ni x}} \Big(\dashint_{C}
 \dashint_{C}|f(z)-f(y)|^{\gamma}\mu(dz)\mu(dy)\Big)^{1/\gamma}.
$$

\begin{theorem}[Fefferman-Stein]
                                            \label{theorem 12.19.1}
 
Let $p\in(1,\infty)$ and $w\in A_{p}$. Then      
for any $f\in \cL_{p}(w)$ we have 
\begin{equation}
                                      \label{12.19.3}
\|f\|_{\cL_{p}(w)}
\leq N\|f^{\#}_{\gamma}\|_{\cL_{p}(w)},
\end{equation}
where  
the constant $N$ depends only on $p$, $N_{0}$, $\gamma$, and
$[w]_{A_{p}}$.
\end{theorem}

Proof. 
In Lemma \ref{lemma 12.16.1} we take $u=f$, $v=u^{C}=|f|$, and
$g=f^{\#}_{\gamma}$. Then we plug $\lambda^{1/p}$  
in place of $\lambda$ in \eqref{6.29.1} and integrate with respect to $\lambda
\in(0,\infty)$.
This yields
$$
\|f\|_{\cL_{p}(w)}^{p}\leq
N\int_{\Omega}(f^{\#}_{\gamma})^{\gamma\beta}(\cM f)^{p-\gamma\beta}\,w(dx).
$$
By using H\"older's inequality, we obtain
$$
\|f\|_{\cL_{p}(w)}^{p}
\leq N\|f^{\#}_{\gamma}\|_{\cL_{p}(w)}^{\gamma\beta}
\|\cM f\|^{p-\gamma\beta}_{\cL_{p}(w)}.
$$
After that 
  it only remains to use
 \eqref{12.19.4}. The theorem is proved.
\begin{remark}
                   \label{remark 12.30.1}
As is easy to see the condition $w\in A_{p}$
can be replaced by $w\in A_{\infty}$.
This is, actually, the form in which
Theorem \ref{theorem 12.19.1} is proved in
\cite{DKr_18}.

\end{remark}

\mysection{An application to fully nonlinear 
parabolic equations}
                                      \label{section 12.20.1}

Here we consider functions on
$$
\bR^{d+1}_{+,+}=\{(t,x)=(t,x_{1},x'):t\geq 0,
x_{1}\geq0, x'\in\bR^{d-1}\}.
$$
We concentrate on parabolic equations in
  $\bR^{d+1}_{+,+}$
 with
 zero Dirichlet
boundary condition.

Denote by $\bS$ the set of symmetric $d\times d$ matrix,
fix $\delta\in(0,1)$, 
and by $S_{\delta}$ denote the subset of $\bS$
consisting of matrices whose eigenvalues are between $\delta$
and $\delta^{-1}$.

Let $A$ be a countable set and suppose that
on $\bR^{d+1}_{+,+}$ for each $\alpha\in A$
we are given an $S_{\delta}$-valued measurable function  
$a^{\alpha}(t,x)
=(a^{\alpha}_{ij}(t,x))$. For $\sfu''=(\sfu''_{ij})\in\bS$
and $(t,x)\in \bR^{d+1}_{+,+}$ introduce
$$
F(\sfu'',t,x)=\sup_{\alpha\in A}\sum_{i,j=1}^{d}
a^{\alpha}_{ij}(t,x)\sfu''_{ij}.
$$
For functions $u=u(t,x)$ having two derivatives in $x$
set $F[u]=F[u](t,x)=F(D^{2}u(t,x),t,x)$.
Also denote
$\partial_{t}=\partial/\partial t$.

For $(t,x)\in \bR^{d+1}_{+,+}$ and $r>0$ denote $B_{r}=
\{x\in\bR^{d}:|x|<r\}$, $B_r(x)=x+B_{r}$,
$B_r^+(x)=B_r(x)\cap\{x^{1}>0\}$,
$$
C^+_{r}(t,x)=[t,t+r^{2})\times B^+_{r}(x) .
$$

Next assumption contains a parameter $\theta\in(0,1)$
which will be specified later.
\begin{assumption}
                              \label{assumption 12.20.1}
There is an $R_{0}\in(0,\infty)$ such that
for any $\alpha\in A$,  $z\in\bR^{d+1}_{+,+} $,
 and $r\in(0,R_{0}]$,   one
can find
$\bar{a}^{\alpha}\in\bS_{\delta}$ (independent of $(t,x)$)
such that
$$
 \dashint_{C^{+}_{r}(z)}\sup_{\alpha\in A}
\big|a^{\alpha}(t,x)-\bar{a}^{\alpha}\big|\,dxdt\leq\theta.
 $$

\end{assumption}


Here is one of the results from \cite{DKr_18}
obtained by combining the above results
and results from \cite{Kr_18}. 

\begin{theorem}
                                   \label{theorem 5.1.1}
 Let $ p_{1},p_{2},p_{3}>d+1$, and
 $u\in W^{ 1, 2}_{1,\loc}(\bR^{d+1}_{+,+})$.
Suppose that $u$ vanishes  on $\{x_{1}=0\}$.
Finally,
take $q\in(-1,p_{1} /(d+1) -1)$. 
 Then there exists
  $\theta=\theta(d,\delta, p_{1},p_{2},p_{3},  q)
\in(0,1) $
such that, if Assumption  \ref{assumption 12.20.1}
is satisfied with this $\theta$, then
$$
\int_{0}^{\infty}\Big(\int_{\bR^{d-1}}\Big(\int_{0}^{\infty}
 x_{1}^{q}\big[|D^{2}u|+|D u| +| u|
\big]^{p_{1}}\,dx_{1}\Big)^{p_{2}/p_{1}}
dx'\Big)^{p_{3}/p_{2}}dt
$$
\begin{equation}
                                               \label{5.1.1}
\le N\int_{0}^{\infty}\Big(\int_{\bR^{d-1}}\Big(\int_{0}^{\infty}
 x_{1}^{q}|\partial_{t}u+F[u]-u|^{p_{1}}\,dx_{1}\Big)^{p_{2}/p_{1}}
dx'\Big)^{p_{3}/p_{2}}dt,
\end{equation}
provided that the left-hand side is finite,
where $N$ depends only on
$d$, $\delta$, $ p_{1}$, $p_{2}$, $p_{3}$, $  q$,  and $R_{0}$.

\end{theorem}

The one-dimensional
 example of $F[u]=D^{2}u$ and $u(t,x)=\sinh x$ shows that \eqref{5.1.1}
is wrong without the additional assumption on its
left-hand side.

\begin{remark}
The reader understands that one has similar estimates
for the integrals with respect to $x_{1}$, $x'$, and $t$
mixed in any other order (cf.~Remark \ref{remark 12.21.1}). For some readers
$\partial_{t}u+F[u]-u$ may look unusual
in comparison with 
$F[u]-\partial_{t}u-u$. One gets the corresponding result for the latter operator
by changing $t\to-t$.
By moving the origin on the $t$-axis to $-\infty$ the reader
obtains the corresponding result when
the integration in \eqref{5.1.1}
with respect to $t$ is extended
to $\bR$. On the other hand, from \eqref{5.1.1} as is
one obtains the corresponding estimates 
for the Cauchy problem as in Theorem 2.5.3 of \cite{Kr_08}.
 
\end{remark}

\begin{remark}
In \cite{KN_09} the authors consider linear $F$ with  coefficients
depending only on time in a measurable way and prove
a priori estimates similar to the one in Theorem \ref{theorem 5.1.1},
however, for any $p_{1}=p_{2}, p_{3}>1$ and $q\in(-1,2p_{1}-1)$.
 The latter range
is much wider than ours $(-1,p_{1} /(d+1) -1)$,
 but our operators are much more general
and we have three   integrals.

Also note that the range $(p_{1}-1,2p_{1}-1)$ was used
in \cite{Kr_01_1} to   build the solvability theory
of parabolic equations  in Sobolev spaces with weights
with the highest order of derivatives being
an arbitrary given number:
 positive, negative,
integral or fractional.

\end{remark}

{\bf Acknowledgement}. The author
is sincerely grateful to Hongjie Dong who read  the first draft of the paper and whose comments
allowed  him  to avoid several errors and
misrepresentations.
  The author also wish to thank P. Stinga
and J. Sauer for
fruitful discussion and providing additional relevant references.

\end{document}